\documentclass[12pt]{article}

\usepackage{amsmath,amsthm,amsfonts,amssymb, color}

\newtheorem{theorem}{Theorem}[section]
\newtheorem{corollary}[theorem]{Corollary}
\newtheorem{example}[theorem]{Example}

\newtheorem{lemma}[theorem]{Lemma}

\newtheorem{remark}[theorem]{Remark}

\def\cD{\mathcal{D}}
\def\cE{\mathcal{E}}

\def\bR{\mathbb{R}}

\def\e{{\varepsilon}}
\def\es{{\varepsilon}^{*}}

\def\lM{\lambda_{\rm max}}
\def\lm{\lambda_{\rm min}}

\def\a{{\bf a}}
\def\b{{\bf b}}
\def\g{{\bf g}}
\def\v{{\bf v}}
\def\w{{\bf w}}
\def\x{{\bf x}}
\def\y{{\bf y}}

\def\A{{\bf A}}
\def\As{{\bf A}^{*}}
\def\Ab{\overline{\bf A}}
\def\B{{\bf B}}
\def\D{{\bf D}}
\def\F{{\bf F}}
\def\G{{\bf G}}
\def\H{{\bf H}}
\def\Hs{{\bf H}^{*}}
\def\I{{\bf I}}
\def\M{{\bf M}}

\def\R{{\bf R}}
\def\Ri{{\bf R}_i(\alpha)}
\def\Rs{{\bf R}^{*}}
\def\T{{\bf T}}
\def\U{{\bf U}}
\def\V{{\bf V}}
\def\Vs{{\bf V}^{*}}
\def\X{{\bf X}}
\def\Y{{\bf Y}}

\begin{document}

\title{Asymptotic theory for longitudinal data with missing responses adjusted by inverse probability weights}

\author{Raluca M. Balan\footnote{Corresponding author. Department of Mathematics and Statistics, University of Ottawa,
585 King Edward Avenue, Ottawa, ON, K1N 6N5, Canada. E-mail
address: rbalan@uottawa.ca} \footnote{Research supported by a
grant from the Natural Sciences and Engineering Research Council
of Canada.}\and
Dina Jankovic\footnote{Department of Mathematics and Statistics, University of Ottawa,
585 King Edward Avenue, Ottawa, ON, K1N 6N5, Canada. E-mail address: djank090@uottawa.ca}}

\date{March 15, 2018}
\maketitle

\begin{abstract}
\noindent In this article, we propose a new method for analyzing longitudinal data which contain responses that are missing at random. This method consists in solving the generalized estimating equation (GEE) of \cite{liang-zeger86} in which the incomplete responses are replaced by values adjusted using the inverse probability weights proposed in \cite{yi-ma-carroll12}. We show that the root estimator is consistent and asymptotically normal, essentially under the some conditions on the marginal distribution and the surrogate correlation matrix as those presented in \cite{XY} in the case of complete data, and under minimal assumptions on the missingness probabilities. This method is applied to a real-life dataset taken from \cite{sommer-katz-tarwotjo84}, which examines the incidence of respiratory disease in a sample of 250 pre-school age Indonesian children which were examined every 3 months for 18 months, using as covariates the age, gender, and vitamin A deficiency.
\end{abstract}

\noindent {\em MSC 2010:} Primary 62F12; Secondary 62J12


\vspace{1mm}

\noindent {\em Keywords:} longitudinal data, generalized estimating equations, asymptotic properties, missing at random, inverse probability weights

\section{Introduction}

Longitudinal data sets are encountered frequently in biostatistics when repeated measurements are made on the same individual. Due to their complexity, the analysis of such data sets presents many challenges for statisticians. Often, one is interested to analyze the relationship between a response variable (for instance the presence of lung cancer) and several explanatory variables (for instance the age, smoking status or family income). In this case, a commonly used method (introduced by Liang and Zeger in the seminal article \cite{liang-zeger86}) is to assume that the marginal distribution of each response follows a generalized linear model (GLM) with regression parameter $\beta$, while the correlation between the responses is modeled by a surrogate correlation matrix which depends on another parameter $\alpha$. The goal of this method is to obtain a consistent estimator of $\beta$, defined as the root of the {\em generalized estimating equation} (GEE). We refer the reader to \cite{mccullagh-nelder89} for a comprehensive account on GLMs, and to \cite{DLZ} for more details about longitudinal data.

Building upon earlier work of \cite{chen-hu-ying99} and \cite{yuan-jennrich98} for estimating equations for classical datasets, the article \cite{XY} contains a thorough analysis of the asymptotic properties of the GEE estimator, including the case when the number of observations made on each individual (called the cluster size) goes to infinity. Similar theoretical  investigations were pursed in \cite{balan-schiopu05} for fixed cluster size, for an estimator defined as the root of a pseudo-likelihood equation, which contains an estimator of the correlation matrix based on the data.

In the presence of incomplete observations, the analysis of longitudinal data becomes even more complex.
Several methods for dealing with longitudinal data which contain missing responses (or missing covariates, or both) have been proposed by various authors. We refer the reader to \cite{chen-yi-cook10, wang-et-al, yi-liu-wu11, yi-ma-carroll12, yi-tan-li15} for a sample of relevant references.

The goal of the present article is to adapt the GEE method of \cite{liang-zeger86, XY} to the case when the responses are {\em missing at random} (a term whose meaning will be explained below). For this, we will replace the incomplete responses by values adjusted using the inverse probability weights proposed in \cite{yi-ma-carroll12}. Under minimal assumptions on the missingness probabilities, we will show that the root estimator of $\beta$ is consistent and asymptotically normal, under essentially the some conditions on the marginal distribution and the surrogate correlation matrix as in \cite{XY}.

We say few words about the notation. We use the convention of omitting the true parameter $\beta_0$ when it is the argument as a function. For instance, we write $\e_i$ instead of $\e_i(\beta_0)$. For sequences $(X_n)_{n \geq 1}$ and $(Y_n)_{n \geq 1}$ of random variables with $Y_n \not=0$ for all $n \geq 1$, we write $X_n=O_p(Y_n)$ if the sequence $(X_n/Y_n)_{n\geq 1}$ is bounded in probability, i.e. for any $\e>0$ there exists $M_{\e}>0$ and an integer $N_{\e}\geq 1$ such that $P(|X_n/Y_n|\leq M_{\e})>1-\e$ for all $n \geq N_{\e}$. We write $X_n=o_{p}(Y_n)$ if $X_n/Y_n \stackrel{p}{\to}0$, where $\stackrel{p}{\to}$ denotes convergence in probability. We write $X_n \stackrel{d}{\to} X$ if $(X_n)_{n\geq 1}$ converges in distribution to $X$.

We conclude the introduction by recalling some basic facts about matrix analysis. We refer the reader to \cite{schott97} for more details.
We denote by ${\rm diag}(\v)$ the diagonal matrix with entries given by $\v=(v_1,\ldots,v_m)$. We denote by $\|\x\|$ the Euclidean norm of a vector $\x$.
If $\A$ is a symmetric matrix, then all its eigenvalues are real. In this case, we write $\A \geq 0$ if $\x^T \A \x \geq 0$ for any vector $\x$, and $\A>0$ if $\x^T \A \x>0$ for any vector $\x$.
For a symmetric $p\times p$ matrix $\A$, we use the following inequality: (see Theorem 3.15 of \cite{schott97})
\begin{equation}
\label{ineq}
\lm(\A) \x^T \x \leq \x^T \A \x \leq \lM(\A)  \x^T \x \quad \mbox{for any} \quad \x \in \bR^p,
\end{equation}
where $\lm(\A)$ is the minimum eigenvalue of $A$ and $\lM(\A)$ is the maximum eigenvalue of $\A$. In particular, $\A\geq 0$ if and only if $\lm(\A)\geq 0$, and $\A>0$ if and only if $\lm(\A)>0$.
If $\A\geq 0$, then $\lM(\A^2)=[\lM(\A)]^2$. If $\A>0$, $\lM(\A^{-1})=1/\lm(\A)$. The left square-root of a matrix $\A>0$ is the matrix $\A^{1/2}$ such that $\A^{1/2}(\A^{1/2})^T=\A$. We let $\A^{-1/2}=(\A^{1/2})^{-1}$.
We let $\|\A\|=\sup_{\|x\|=1}\|\A \x\|=\big(\lM(\A^T \A)\big)^{1/2}$ be the spectral norm of a matrix $\A=(a_{ij})_{i\leq n, j\leq m}$, which is equivalent to its Euclidean norm given by $\|\A\|_{E}=\big(\sum_{i=1}^{n}\sum_{j=1}^{m}a_{ij}^2\big)^{1/2}$. We denote by ${\rm tr}(\A)$ the trace of matrix $\A$ and by ${\rm det}(\A)$ its determinant. Note that ${\rm tr}(\A)$ is the sum of the eigenvalues of $\A$ and ${\rm det}(\A)$ is the product of the eigenvalues of $\A$ (see Theorem 3.5 of \cite{schott97}).

This article is organized as follows. In Section \ref{section-EE} we introduce our framework, we define the estimating equation and we discuss some of its properties. In Section \ref{section-normal}, we prove the consistency and asymptotic normality of the root estimator, under essentially the same conditions as in \cite{XY} in the case of complete data. The most complicated of these conditions involves the derivative of the estimating equation and is called condition (CC). In Section \ref{section-CC}, we give some sufficient conditions for (CC).
In Section \ref{section-ex}, we apply our method to a real-life data set. Finally, Appendix A contains some auxiliary results which are used in Section \ref{section-CC}.

\section{The estimating equation}
\label{section-EE}

In this section, we introduce our framework and we define a generalized estimating equation which can be used when some of the responses are missing at random.

We consider $n$ individuals whose measurements are recorded on $m$ occasions. For each $i=1,\dots,n$ and $j=1,\dots,m$, we denote by $Y_{ij}$ the response of individual $i$ at time $j$. Some of these responses are missing. We let
$$I_{ij}=
\begin{cases}
    1,& \mbox{if $Y_{ij}$ is observed}\\
    0,& \mbox{if $Y_{ij}$ is missing}
\end{cases}$$
We let $\textbf{Y}_{i}=(Y_{i1},\ldots,Y_{im})^T$ be the vector of responses of the $i^{th}$ individual and $\textbf{I}_{i}=(I_{i1},\ldots,I_{im})^T$ be the vector of missingness indicators for this individual.

For each $i=1,\dots,n$, and $j=1,\dots,m$, we let
$\X_{ij}=(X^{(1)}_{ij},\ldots,X^{(p)}_{ij})^T$
be the $p$-dimensional vector of covariates for individual $i$ at time $j$. We assume that $\X_{ij}$ is random. The following $m\times p$ matrix contains the covariates of the $i^{th}$ individual:
$$
\textbf{X}_i =
  \begin{bmatrix}
    \X^T_{i1} \\
    \vdots \\
    \X^T_{im}
  \end{bmatrix}
  = \begin{bmatrix}
  X^{(1)}_{i1} & \hdots & X^{(p)}_{i1}\\
  \vdots & \ddots & \vdots \\
  X^{(1)}_{im} & \hdots & X^{(p)}_{im}
  \end{bmatrix}
$$

We assume that ${\{(\textbf{Y}_{i},\textbf{X}_i,\textbf{I}_i)\}}_{i\ge1}$ are independent and identically distributed (i.i.d.), and there exists a one-to-one differentiable function $\mu$ on $\mathbb{R}$ such that
\begin{equation}
\label{model-ass}
\mu_{ij}(\beta):=E(Y_{ij}|\textbf{X}_i)=\mu(X^T_{ij}\beta) \quad \mbox{and} \quad \sigma_{ij}^2(\beta):={\rm Var}(Y_{ij}|\textbf{X}_i)=\phi\mu'(X^T_{ij}\beta),
\end{equation}
for a $p$-dimensional parameter $\beta$ and a nuisance parameter $\phi$. In the present article, we will assume that $\phi=1$.
The inverse $g$ of the function $\mu$ is called the {\em link function}. Let $\mu_i(\beta)=(\mu_{i1}(\beta),\dots,\mu_{im}(\beta)^T$. We denote by $\D_i(\beta)$ the $m\times p$ matrix:
$$\textbf{D}_i(\beta)=\frac{\partial\mu_i(\beta)}{\partial\beta^T}=
  \begin{bmatrix}
    \frac{\partial\mu_{i1}}{\partial\beta^T}(\beta) \\
    \vdots \\
    \frac{\partial\mu_{im}}{\partial\beta^T}(\beta)
  \end{bmatrix}$$

\noindent Note that $\frac{\partial\mu_{ij}}{\partial\beta^T}(\beta)=\textbf{X}^T_{ij}\mu'(\textbf{X}^T_{ij}\beta)$ and hence $\textbf{D}_i(\beta)=\textbf{A}_i(\beta)\textbf{X}_i$, where $\textbf{A}_i(\beta)$ is the diagonal matrix with entries $\mu'(\textbf{X}^T_{i1}\beta),\dots,\mu'(\textbf{X}^T_{im}\beta)$ for $j=1,\ldots,m$.

\begin{example}{\rm (Normal Linear Regression for Quantitative Responses)
When the responses $Y_{ij}$ represent quantitative measurements, we may assume that $Y_{ij}$ has an normal distribution with mean $\mu_{ij}(\beta)=\X_{ij}^T\beta$ and known variance $\sigma_{ij}^2(\beta)=\phi$, for a nuissance parameter $\phi$ which is estimated separately. In this case, $\mu(x)=x$ and $\mu'(x)=1$. The link function is $g(x)=x$.}
\end{example}

\begin{example}{\rm (Log-linear Regression for Count-type Responses)
When the responses $Y_{ij}$ represent count-type measurements, we may assume that $Y_{ij}$ has a Poisson distribution with mean $\mu_{ij}(\beta)=\exp(\X_{ij}^T\beta)$. In this case, $\mu(x)=e^x$ and $\mu'(x)=e^x$. The link function is $g(x)=\log x$ for $x>0$.}
\end{example}

\begin{example}\label{example3}{\rm (Logistic Regression for Binary Responses)
When the responses $Y_{ij}$ represent binary measurements, we may assume that $Y_{ij}$ has a Bernoulli distribution with mean $\mu_{ij}(\beta)=\frac{\exp(\X_{ij}^T\beta)}{1+\exp(\X_{ij}^T\beta)}$. In this case, $\mu(x)=\frac{e^x}{1+e^x}$ and $\mu'(x)=\frac{e^x}{(1+e^x)^2}$. The link function is $g(x)=\log \frac{x}{1-x}=:{\rm logit}(x)$ for $x \in (0,1)$.}
\end{example}

We consider the following marginal model:
$$Y_{ij}=\mu_{ij}(\beta)+\varepsilon_{ij}(\beta) \quad j=1, \dots ,m.$$

We let $\e_i(\beta)=(\e_{i1}(\beta),\ldots,\e_{im}(\beta))^T$ be the residuals, for $i=1,\ldots,n$. Let $\Sigma_i(\beta)=(\sigma_{i,jk}(\beta))_{1\leq j,k\leq m}$ be the conditional covariance matrix of $\textbf{Y}_i$ given $\textbf{X}_i$, with entries $\sigma_{i,jk}(\beta)=\mathrm{E}[\varepsilon_{ij}(\beta)\varepsilon_{ik}(\beta)|\textbf{X}_i]$.
In particular, $\sigma_{i,jj}(\beta)=\sigma^2_{ij}(\beta)=\mu'(X^T_{ij}\beta)$.
Note that
$$\Sigma_i(\beta)=\textbf{A}_i(\beta)^{1/2}\overline{\textbf{R}}_i\textbf{A}_i(\beta)^{1/2},$$
where $\overline{\textbf{R}}_i=(\overline{r}_{i,jk})_{1\le j,k \le m}$ is the conditional correlation matrix of $\textbf{Y}_i$ given $\textbf{X}_i$ with entries:
$\overline{r}_{i,jk}=\sigma_{i,jk}(\beta)/[\sigma_{ij}(\beta)\sigma_{ik}(\beta)]$.


The method proposed in \cite{liang-zeger86} consists in replacing the unknown correlation matrix $\overline{\textbf{R}}_i$ by a surrogate correlation matrix $\textbf{R}_i(\alpha)$ depending on a parameter $\alpha$ (to be estimated separately), and solving the Generalized Estimated Equation (GEE):

\begin{equation} \label{eq1}
\displaystyle\sum_{i=1}^n\textbf{D}^T_i(\beta)\textbf{V}^{-1}_i(\beta,\alpha)(\textbf{Y}_i-\mu_i(\beta))=0,
\end{equation}

\noindent where $\textbf{V}_i(\beta)=\textbf{A}_i(\beta)^{1/2}\textbf{R}_i(\alpha)\textbf{A}_i(\beta)^{1/2}$ is an approximation of the unknown covariance matrix $\Sigma_i(\beta)$. Equation \eqref{eq1} can be written equivalently as:
\begin{equation}
\label{new-eq1}
\sum_{i=1}^n \X_i^T \A_i(\beta)^{1/2} \Ri^{-1}\A_i(\beta)^{-1/2}(\textbf{Y}_i-\mu_i(\beta))=0.
\end{equation}
Typical examples of matrices $\Ri$ are:
$$a) \ \textbf{R}_i(\alpha)=
\begin{bmatrix}
1        & \alpha_1 &        0 & 0      & \hdots & 0 & 0\\
\alpha_1 &        1 & \alpha_2 & 0      & \hdots & 0 & 0\\
\vdots   & \vdots   & \ddots   & \vdots & \vdots & 0 & 0 \\
0        &        0 &        0 & 0      & \hdots & 1 & \alpha_{m-1} \\
0        &        0 &        0 & 0      & \hdots & \alpha_{m-1} & 1
\end{bmatrix}
\quad \mbox{or} \quad
b) \ \textbf{R}_i(\alpha)=
\begin{bmatrix}
1 & \alpha & \hdots & \alpha \\
\alpha & 1 & \hdots & \alpha \\
\vdots & \vdots & \ddots & \vdots \\
\alpha & \alpha & \hdots & 1
\end{bmatrix}
$$
Case {\em a)} is called {\em $1$-dependent}, whereas case {\em b)} is called {\em exchangeable}.

The case $\textbf{R}_i(\alpha)=\textbf{I}$ for all $i=1,\dots,n$ is called {\em working independence}. In this case, equation (\ref{new-eq1}) becomes
\begin{equation}\label{twooneoneprime}
\displaystyle\sum_{i=1}^{n}\displaystyle\sum_{j=1}^{m}\X_{ij}(Y_{ij}-\mu(\X_{ij}^T\beta))=0.
\end{equation}

In 2003, Xie and Yang proved rigorously in \cite{XY} that equation (\ref{eq1}) has a root $\widehat{\beta}_n$ which is a consistent estimator of $\beta$, and derived the asymptotic normality of this estimator.
In this article, we develop a method similar to that of \cite{XY} which can be applied when some of the responses are missing.

We assume that the the responses are {\em missing at random} (MAR), i.e.
\begin{center}
$\Y_i$ and $\I_i$ are conditionally independent given $\X_i$, for any $i=1,\ldots,n$.
\end{center}

For any $i=1,\dots,n$ and $j=1,\dots,m$, we let
$\pi_{ij}=P(I_{ij}=1|\textbf{X}_i,\textbf{Y}_i)=P(I_{ij}=1|\textbf{X}_i)$.
Then $E(I_{ij}|\textbf{X}_i)=E(I^2_{ij}|\textbf{X}_i)=\pi_{ij}$.
We consider the {\em inverse probability weighted response}
\begin{equation}
\label{def-Ystar}
Y^*_{ij}=\frac{Y_{ij}I_{ij}}{\pi_{ij}}.
\end{equation}

We let $\textbf{Y}^*_i=(Y^*_{i1},\dots,Y^*_{im})^T$ be the vector of weighted responses for the $i^{th}$ individual and $\varepsilon^*_i(\beta)=(\varepsilon^*_{i1}(\beta),\dots,\varepsilon^*_{im}(\beta))^T$, where $\varepsilon^*_{ij}(\beta)=Y^*_{ij}-\mu_{ij}(\beta)$, $j=1,\dots,m$.

The next result shows that the weighted response $Y_{ij}^*$ has the same mean as the original response $Y_{ij}$.

\begin{lemma} \label{lem1}
For each $i=1,\dots,n$ and $j=1,\dots,m$, $\mathrm{E}(Y^*_{ij}|\X_i)=\mu_{ij}(\beta)$. Therefore $$\mathrm{E}(\varepsilon^*_i(\beta)|\textbf{X}_i)=0 \ for \ all \ i \ge 1.$$
\end{lemma}

\begin{proof}
Note that
\begin{equation}
\label{eq2}
Y^*_{ij}-Y_{ij}=\Big(\frac{I_{ij}}{\pi_{ij}}-1\Big)Y_{ij}.
\end{equation}
Using (\ref{eq2}) and double conditioning, we have:
\begin{align*}
\mathrm{E}(Y^*_{ij}|\textbf{X}_i)&=
\mathrm{E}(Y^*_{ij}-Y_{ij}|\textbf{X}_i)+\mathrm{E}(Y_{ij}|\textbf{X}_i)
=\mathrm{E}\Big[\Big(\frac{I_{ij}}{\pi_{ij}}-1\Big)Y_{ij}|\textbf{X}_i\Big]+\mu_{ij}(\beta)\\
&=\mathrm{E}\Big[Y_{ij}\mathrm{E}\Big[\frac{I_{ij}}{\pi_{ij}}-1|\textbf{X}_i,\textbf{Y}_i\Big]|
\textbf{X}_i\Big]+\mu_{ij}(\beta)=\mu_{ij}(\beta),
\end{align*}
\noindent where for the last line, we used the fact that,
\begin{equation}\label{eq3}
\mathrm{E}\Big[\frac{I_{ij}}{\pi_{ij}}-1|\textbf{Y}_i,\textbf{X}_i\Big]=\mathrm{E}\Big[\frac{I_{ij}}{\pi_{ij}}-1|\textbf{X}_i\Big]=\frac{1}{\pi_{ij}}\mathrm{E}(I_{ij}|\textbf{X}_i)-1=0.
\end{equation}
due to the (MAR) assumption.
\end{proof}

For each $i=1,\dots,n$ and $j,k=1,\dots,m$, we consider the probability that both responses $Y_{ij}$ and $Y_{ik}$ are missing, given $\X_i$:
\begin{center}
$q_{i,jk}=P(I_{ij}=1, I_{ik}=1|\textbf{X}_i).$
\end{center}

In the next lemma, we compute the conditional covariance matrix of $\Y^*_i$ given $\X_i$. Note that the expected value of a matrix $\textbf{A}=(A_{jk})_{1 \le j,k \le m}$ whose elements are random variables $A_{jk}$ is, by definition, the matrix $\mathrm{E}(\textbf{A})=\{{\mathrm{E}(A_{jk})}\}_{1 \le j,k \le m}$.

\begin{lemma}
The conditional covariance matrix of $\Y^*_i$ given $\textbf{X}_i$ is
$$\Sigma^*_i(\beta):=\mathrm{E}[\varepsilon^*_i(\beta)\varepsilon^*_i(\beta)^T|\textbf{X}_i]=(\sigma^*_{i,jk}(\beta))_{1\le j,k \le m},$$
where
\begin{equation}
\label{def-cov-star}\sigma_{i,jk}^*(\beta)=\sigma_{i,jk}(\beta)+\left(\frac{q_{i,jk}}{\pi_{ij}\pi_{ik}}-1 \right)\big(\sigma_{i,jk}(\beta)+\mu_{ij}(\beta)\mu_{ik}(\beta)\big).
\end{equation}
In particular, for any $j=1,\ldots,m$, the conditional marginal variance of $Y_{ij}^*$ given $\X_i$ is:
\begin{align}
\nonumber
\sigma_{i,jj}^*(\beta)&=\sigma_{ij}^2(\beta)+\left(\frac{1}{\pi_{ij}}-1 \right)(\sigma_{ij}^2(\beta)+\mu_{ij}^2(\beta))\\
\label{sigma_star_ijj}
&=\mu'(\X_{ij}^T \beta)+\left(\frac{1}{\pi_{ij}}-1 \right)(\mu'(\X_{ij}^T\beta)+\mu^2(\X_{ij}^T \beta)).
\end{align}
\end{lemma}

\begin{proof}
For any $j,k=1,\dots,m$ fixed,
\begin{align}
\nonumber\sigma^*_{i,jk}(\beta)\nonumber&
=\mathrm{E}[(Y^*_{ij}-\mu_{ij}(\beta))(Y^*_{ik}-\mu_{ik}(\beta))|\textbf{X}_i]
\\
\nonumber
\nonumber
&=\mathrm{E}[(Y^*_{ij}-Y_{ij})(Y^*_{ik}-Y_{ik})|\textbf{X}_i]+
\mathrm{E}[(Y^*_{ij}-Y_{ij})(Y_{ik}-\mu_{ik}(\beta))|\textbf{X}_i]\\
&+\mathrm{E}[(Y^*_{ik}-Y_{ik})(Y_{ij}-\mu_{ij}(\beta))|\textbf{X}_i]+
\mathrm{E}[(Y_{ij}-\mu_{ij}(\beta))(Y_{ik}-\mu_{ik}(\beta))|\textbf{X}_i].\label{eq4}
\end{align}

\noindent We treat separately the four terms. By \eqref{eq2}, the second term is equal to
\begin{align*}
\mathrm{E}[(Y^*_{ij}-Y_{ij})(Y_{ik}-\mu_{ik}(\beta))|\textbf{X}_i]
=\mathrm{E}\Big[Y_{ij}(Y_{ik}-\mu_{ij}(\beta))\mathrm{E}\Big[\frac{I_{ij}}{\pi_{ij}}-1|
\textbf{X}_i,\textbf{Y}_i\Big]|\textbf{X}_i\Big]=0,
\end{align*}
using (\ref{eq3}) for the last equality. Similarly, the third term in (\ref{eq4}) is also equal to 0. Note that the fourth term in (\ref{eq4}) is equal to
$\sigma_{i,jk}(\beta)$. Hence,
\begin{equation}\label{eq5}
\sigma^*_{i,jk}(\beta)=\mathrm{E}[(Y^*_{ij}-Y_{ij})(Y^*_{ik}-Y_{ik})|\textbf{X}_i]+\sigma_{i,jk}(\beta).
\end{equation}

\noindent By \eqref{eq2},
\begin{equation}
\label{eq6}
\mathrm{E}[(Y^*_{ij}-Y_{ij})(Y^*_{ik}-Y_{ik})|\textbf{X}_i]
=\mathrm{E}\Big[Y_{ij}Y_{ik}\mathrm{E}\Big[\Big(\frac{I_{ij}}
{\pi_{ij}}-1\Big)\Big(\frac{I_{ik}}{\pi_{ik}}-1\Big)|\textbf{X}_i,\textbf{Y}_i\Big]|
\textbf{X}_i\Big].
\end{equation}

\noindent We compute separately the inner conditional expectation. By the (MAR) assumption,
\begin{align*}
\mathrm{E}\Big[\Big(\frac{I_{ij}}{\pi_{ij}}-1\Big)
\Big(\frac{I_{ik}}{\pi_{ik}}-1\Big)|\textbf{X}_i,\textbf{Y}_i\Big]&=
\quad \frac{1}{\pi_{ij}\pi_{ik}}\mathrm{E}(I_{ij}I_{jk}|\textbf{X}_i) - \frac{1}{\pi_{ij}}\mathrm{E}(I_{ij}|\textbf{X}_i) - \frac{1}{\pi_{ik}}\mathrm{E}(I_{ik}|\textbf{X}_i) + 1\\
&=\frac{1}{\pi_{ij} \pi_{ik}}q_{i,jk}-1.
\end{align*}

\noindent Coming back to (\ref{eq6}), we obtain
\begin{align}
\nonumber\mathrm{E}[(Y^*_{ij}-Y_{ij})(Y^*_{ik}-Y_{ik})|\textbf{X}_i]&=
\Big(\frac{1}{\pi_{ij}\pi_{ik}}q_{i,jk}-1\Big)\mathrm{E}[Y_{ij}Y_{ik}|\textbf{X}_i]\\
&=\Big(\frac{1}{\pi_{ij}\pi_{ik}}q_{i,jk}-1\Big)\big(\sigma_{i,jk}+\mu_{ij}(\beta)\mu_{ik}(\beta)\big).
\label{product}
\end{align}
Relation \eqref{def-cov-star} follows from relations (\ref{eq5}) and (\ref{product}).
The last statement follows from \eqref{def-cov-star} and our model assumptions \eqref{model-ass}, using the fact that  $q_{i,jj}=P(I_{ij}=1|\textbf{X}_i)=\pi_{ij}$.
\end{proof}

\vspace{3mm}

Let $\As_i(\beta)$ be the diagonal matrix with entries $\sigma_{i,jj}^*(\beta),j=1,\ldots,m$. Then
\begin{equation}
\label{sigma-R-star}
\Sigma_i^*(\beta)=\As_i(\beta)^{1/2}\Rs_i \As_i(\beta)^{1/2},
\end{equation}
where $\Rs_i=(r_{i,jk}^*)_{1\leq j,k\leq m}$ is the conditional correlation matrix of $\Y_i^*$ given $\X_i$:
$$r_{i,jk}^*=\frac{{\rm Cov}(Y_{ij}^*, Y_{ik}^*|\X_i)}{\sqrt{{\rm Var}(Y_{ij}^*|\X_i)}\cdot\sqrt{{\rm Var}(Y_{ik}^*|\X_i)}}=\frac{\sigma_{i,jk}^*(\beta)}{\sqrt{\sigma_{i,jj}^*(\beta)}\cdot \sqrt{\sigma_{i,kk}^*(\beta)}}.$$

In practice, the matrix $\Rs_i$ is unknown. Following the same idea as in \cite{liang-zeger86} in the case of complete data, we replaced the matrix $\Rs_i$ by a surrogate matrix $\Ri$ which depends on an unknown parameter $\alpha$ (to be estimated separately).
We define
\begin{equation}
\label{def_Vi_star}
\textbf{V}_i^*(\beta,\alpha)=\textbf{A}_i^*(\beta)^{1/2}\textbf{R}_i(\alpha)\textbf{A}_i^*(\beta)^{1/2}.
\end{equation}

We are interested in solving the equation
\begin{equation}
\label{eq7}
\g_n(\beta):=\displaystyle\sum_{i=1}^n\textbf{D}_i^T(\beta)\textbf{V}_i^*(\beta,\alpha)^{-1}(\textbf{Y}_i^*-\mu_i(\beta))=0.
\end{equation}
Note that equation (\ref{eq7}) is the analogue of equation \eqref{eq1} for the case of missing responses which are adjusted using the inverse probability weights.
 Note that
\begin{align*}
\g_n(\beta)=\sum_{i=1}^{n}\X_i^T \F_i(\beta)\Ri^{-1}\As_i(\beta)^{-1/2}\es_i(\beta),
\end{align*}
where $\F_i(\beta)=\A_i(\beta)\As_i(\beta)^{-1/2}$ is the diagonal matrix with entries $f_{ij}(\beta),j=1,\ldots,m$:
\begin{equation}
\label{def-fij}
f_{ij}(\beta)=\frac{\sigma_{ij}^2(\beta)}{\sqrt{\sigma_{i,jj}^{*}(\beta)}}=\frac{\mu'(\X_{ij}^T \beta)}{\sqrt{\mu'(\X_{ij}^T \beta)+\left(\frac{1}{\pi_{ij}}-1 \right)\big(\mu'(\X_{ij}^T \beta)+\mu^2(\X_{ij}^T\beta)\big)}}.
\end{equation}

The following result gives the mean and the covariance matrix of $\g_n(\beta)$.

\begin{lemma}
\label{cov-gn-lemma}
$\g_n(\beta)$ is an unbiased estimating function, i.e. $E[\g_n(\beta)]=0$ for all $\beta$.
The covariance matrix of $\g_n(\beta)$ is
$$\M_n(\beta):=E[\g_n(\beta)\g_n(\beta)^T]=E[\M_n^*(\beta)],$$
where
$$\M_n^*(\beta)=\displaystyle\sum_{i=1}^n\D_i^T(\beta)
\Vs_i(\beta,\alpha)^{-1}\Sigma_i^*(\beta)\Vs_i(\beta,\alpha)^{-1}
\D_i(\beta).$$
\end{lemma}

\begin{proof} The first statement follows by Lemma \ref{lem1} since
$$E[\g_n(\beta)]=\sum_{i=1}^{n}E[\D_i^T(\beta)\Vs_i(\beta,\alpha)^{-1} E[\es_i(\beta)|\X_i]]=0.$$

We proceed now with the calculation of the covariance matrix of $\g_n(\beta)$. Note that
\begin{center}
$\g_n(\beta)\g_n(\beta)^T=\displaystyle\sum_{i=1}^n\displaystyle\sum_{l=1}^n\textbf{D}_i^T(\beta)\textbf{V}_i^*(\beta,\alpha)^{-1}\varepsilon_i^*(\beta)\varepsilon_l^*(\beta)^T\textbf{V}_l^*(\beta,\alpha)^{-1}\textbf{D}_l(\beta)$.
\end{center}
Since $\pi_{ij}=P(I_{ij}=1|\textbf{X}_i)=h_j(\textbf{X}_i)$ for a function $h_j$,
$$\varepsilon_{ij}^*(\beta)=Y_{ij}^*-\mu(\X_{ij}^T\beta)=
\frac{Y_{ij}I_{ij}}{\pi_{ij}}-\mu(\X_{ij}^T\beta)=\Phi_j(Y_{ij},\textbf{X}_{i},I_{ij},\beta)$$
for a certain function $\Phi_j$. Since ${\{(\textbf{Y}_i,\textbf{X}_i,\textbf{I}_i)\}}_{i \ge 1}$ are independent, it follows that ${\{\varepsilon_i^*(\beta)\}}_{i \ge 1}$ are independent. The same argument shows that ${\{\textbf{D}_i^T(\beta)\textbf{V}_i^*(\beta,\alpha)^{-1}\varepsilon_i^*(\beta)\}}_{i \ge 1}$ are independent. Note that
\begin{align*}
\mathrm{E}[\textbf{D}_i^T(\beta)\textbf{V}_i^*(\beta,\alpha)^{-1}\varepsilon_i^*(\beta)]&=\mathrm{E}[\mathrm{E}[\textbf{D}_i^T(\beta)\textbf{V}_i^*(\beta,\alpha)^{-1}\varepsilon_i^*(\beta)|\textbf{X}_i]]\\&=\mathrm{E}[\textbf{D}_i^T(\beta)\textbf{V}_i^*(\beta,\alpha)^{-1}\mathrm{E}[\varepsilon_i^*(\beta)|\textbf{X}_i]]\\&=0,
\end{align*}
\noindent where for the last equality we used Lemma \ref{lem1}. Therefore if $i \neq l$,
\begin{center}
$\mathrm{E}[\textbf{D}^T_i(\beta)\textbf{V}^*_i(\beta,\alpha)^{-1}\varepsilon^*_i(\beta)\varepsilon^*_l(\beta)\textbf{V}^*_l(\beta,\alpha)^{-1}\textbf{D}_l(\beta)]=0.$
\end{center}

\noindent Coming back to the calculation of $\mathrm{E}[\g_n(\beta)\g_n(\beta)^T]$, we obtain using conditioning again
\begin{align*}
\mathrm{E}[\g_n(\beta)\g_n(\beta)^T]&=\displaystyle\sum_{i=1}^n\mathrm{E}[\mathrm{E}[\textbf{D}_i^T(\beta)\textbf{V}_i^*(\beta,\alpha)^{-1}\varepsilon_i^*(\beta)\varepsilon_i^*(\beta)^T\textbf{V}_i^*(\beta,\alpha)^{-1}\textbf{D}_i(\beta)|\textbf{X}_i]]\\&=\displaystyle\sum_{i=1}^n\mathrm{E}[\textbf{D}_i^T(\beta)\textbf{V}_i^*(\beta,\alpha)^{-1}\mathrm{E}[\varepsilon_i^*(\beta)\varepsilon_i^*(\beta)^T|\textbf{X}_i]\textbf{V}_i^*(\beta,\alpha)^{-1}\textbf{D}_i(\beta)]\\&=\displaystyle\sum_{i=1}^n\mathrm{E}[\textbf{D}_i^T(\beta)\textbf{V}_i^*(\beta,\alpha)^{-1}\Sigma_i^*(\beta)\textbf{V}_i^*(\beta,\alpha)^{-1}\textbf{D}_i(\beta)]\\&=\mathrm{E}[\textbf{M}_n^*(\beta)].
\end{align*}

\noindent This finishes the proof.
\end{proof}

\begin{remark}\label{remark214}
{\rm Using the fact that $\textbf{D}_i(\beta)=\textbf{A}_i(\beta)\textbf{X}_i$
and relations \eqref{sigma-R-star} and \eqref{def_Vi_star},
we obtain the following alternative formula for $\textbf{M}_n^*(\beta)$:
$$\textbf{M}_n^*(\beta)=\displaystyle\sum_{i=1}^n\textbf{X}_i^T\textbf{F}_i(\beta)\textbf{R}_i(\alpha)^{-1}\textbf{R}_i^*\textbf{R}_i(\alpha)^{-1}
\textbf{F}_i(\beta)\textbf{X}_i.$$

We denote $\tau_n=\max\limits_{i \le n}\lambda_{\rm max}(\textbf{R}_i(\alpha)^{-1/2}\textbf{R}_i^*\textbf{R}_i(\alpha)^{-1/2})=
\max\limits_{i \le n}\lambda_{\rm max}(\textbf{R}_i^{-1}(\alpha)\textbf{R}_i^*)$.
By \eqref{ineq}, $\textbf{M}_n^*(\beta) \le \tau_n\textbf{H}_n^*(\beta)$, where
\begin{equation}
\label{def-Hstar}
\Hs_n(\beta)=\sum_{i=1}^{n}\X_i^T \F_i(\beta) \Ri^{-1}\F_i(\beta)\X_i=\sum_{i=1}^{n}\D_i(\beta)^T \Vs_i(\beta,\alpha)^{-1}\D_i(\beta).
\end{equation}
By taking the expectation on both sides of this inequality, we infer that $$\textbf{M}_n(\beta) \le \tau_n\textbf{H}_n(\beta),$$ where $\textbf{H}_n(\beta)=\mathrm{E}[\textbf{H}_n^*(\beta)]$.
Note that $\tau_n \leq m \widetilde{\lambda}_n$, where $\widetilde{\lambda}_n=\max_{i\leq n}\lM(\Ri^{-1})$.
The advantage of working with $\widetilde{\lambda}_n$ instead of $\tau_n$ is that $\widetilde{\lambda}_n$ does not depend on the unknown correlation matrix $\R_i^{*}$.
}
\end{remark}

\section{Consistency and asymptotic normality}
\label{section-normal}

In this section, we show that under certain conditions, equation (\ref{eq7}) has a solution $\widehat{\beta}_n$ which is a consistent estimator of $\beta$. The proofs are similar to those presented in \cite{XY} in the case of complete data.

We consider the negative derivative of our estimating function $\g_n(\beta)$:
\begin{center}
$\mathcal{D}_n(\beta)=-\frac{\partial}{\partial\beta^T}\g_n(\beta).$
\end{center}
This derivative plays an important form in the present article. Its explicit formula is given in Section \ref{section-CC} below. It is important to note that $\mathcal{D}_n(\beta)$ is non-symmetric. We consider the ball $B_n^*(r)=\{\beta ; \|\textbf{H}_n^{1/2}(\beta - \beta_0)\| \le \tau_n^{1/2}r\}$, where $\beta_0$ is the true value of $\beta$.

\vspace{2mm}

Similarly to \cite{XY}, we consider the following conditions:
\vspace{2mm}

\noindent $(I_w^*)$  $\lm(\H_n)/\tau_n \to \infty$ as $n \to \infty$.\\
$(L_w^*)$ There exists a constant $c_0 > 0$ such that for any $r > 0$,
$$P(\x^T\H_n^{-1/2}\mathcal{D}_n(\beta)\H_n^{-1/2}\x \ge c_0 \ \mbox{for any $\beta \in B_n^*(r)$ and $\x \in \bR^p$ with $\|\x\|=1$}) \to 1.$$
$(D_w^*)$ For any $r>0$,
$P(\mathcal{D}_n(\beta) \ \mbox{is non-singular for any $\beta \in B_n^*(r)$}) \to 1$.

\vspace{3mm}

Under $(I_w^*)$, $\lm(\H_n)>0$ for $n$ large enough. Hence, $\H_n>0$ for $n$ large enough.

The following result shows that under these conditions, there exists an estimator $\widehat{\beta}_n$ which is the root of the equation $\g_n(\beta)=0$ and this estimator is consistent.
\begin{theorem}
\label{existence-th}
Under conditions $(I_w^*)$, $(L_w^*)$ and $(D_w^*)$, we have:\\
a)$ P({\rm there \ exist \ } r>0 \ {\rm and} \ \widehat{\beta}_n \ \in B_n(r) \ {\rm such \ that} \ \g_n(\widehat{\beta}_n)=0) \to 1$\\
b) $\widehat{\beta}_n \stackrel{P}{\to} \beta_0$.
\end{theorem}

\begin{proof}
This use the same argument as in the proof of Theorem 2 of \cite{XY}. We give only the sketch of this argument.

a) Let $\Omega_n^*(r)$ be the event where $\x^T\H_n^{-1/2}\mathcal{D}_n(\beta)\H_n^{-1/2}\x \ge c_0$ for any $\beta \in B_n^*(r)$ and for any $x \in \bR^p$ with $\|\x\|=1$, and $\cD_n(\beta)$ is non-singular for any $\beta \in B_n^*(r)$. By conditions $(L_w^*)$ and $(D_w^*)$, $P(\Omega_n^*(r)) \to 1$ for any $r>0$.

On the event $\Omega_n^*(r)$, the function $\T_n(\beta)=\H_n^{-1/2}\g_n(\beta)$ is one-to-one, since its derivative is non-singular. This function is also differentiable. Let $E_n^*(r)$ be the event where
$\|\T_n(\beta_0)\|\leq \inf_{\beta \in \partial B_n^*(r)}\|\T_n(\beta)-\T_n(\beta_0)\|$, where $\partial B_n^*(r)=\{\beta; \|\H_n^{1/2}(\beta-\beta_0)\|=\tau_n^{1/2}r\}$. Let $\widetilde{\Omega}_n(r)$ be the event that there exists $\widehat{\beta}_n \in B_n^*(r)$ such that $\g_n(\widehat{\beta}_n)=0$. By Lemma A of \cite{chen-hu-ying99},
\begin{equation}
\label{Omega-E}
E_n^*(r) \cap \Omega_n^*(r) \subset \widetilde{\Omega}_n(r).
\end{equation}

Therefore, it suffices to show that for any $\e>0$, there exists $r=r_{\e}>0$ and an integer $N_{\e}\geq 1$ such that
\begin{equation}
\label{eqD}
P(E_n^*(r) \cap \Omega_n^*(r))\geq 1-\e \quad \mbox{for all} \ n \geq N_{\e}.
\end{equation}

Let $\e>0$ be arbitrary and $r=\frac{1}{c_0}\sqrt{\frac{2p}{\varepsilon}}$. By applying Talyor's formula to $\T_n(\beta)$, it can be proved that on the event $\Omega_n^*(r)$,
$\inf_{\beta \in \partial B_n^*(r)}\|\T_n(\beta)-\T_n(\beta_0)\| \geq c_0 \tau_n^{1/2}r$.
Hence, the event $\{\|\T_n(\beta_0)\| \leq c_0 \tau_n^{1/2}r\} \cap \Omega_n^*(r)$ is contained in $E_n^*(r) \cap \Omega_n^*(r)$. By Chebushev's inequality and the choice of $r$, for any $n\geq1$, we have:
$$P(\|\T_n(\beta_0)\| \leq c_0 \tau_n^{1/2}r) \geq 1-\frac{1}{c_0^2\tau_n r^2}E[\|\T_n(\beta_0)\|^2]=1-\frac{\e}{2}.$$
Since $P(\Omega_n^*(r)) \to 1$, there exists an integer $N_{\e}\geq 1$ such that $P(\Omega_n^*(r))>1-\e/2$ for all $n\geq N_{\e}$. Relation \eqref{eqD} follows by Bonferroni's inequality.

b) Let $\delta>0$ and $\e>0$ be arbitrary. Let $r=r_{\e}$ and $N_{\e}$ as in part a). By condition $(I_w^*)$, there exists $N_\delta^* \ge 1$ such that $\frac{\lambda_{\rm min}(\textbf{H}_n)}{\tau_n} \ge \frac{\delta^2}{r}$ for all $n \ge N_\delta^*$.
On the event $\Omega_n^*(r) \cap E_n^*(r)$,
$\|\widehat{\beta}_n - \beta_0\| \le \|\textbf{H}_n^{-1/2}\|\cdot\|\textbf{H}_n^{1/2}(\widehat{\beta}_n - \beta_0)\| \le \Big[\frac{1}{\lambda_{\min} (\textbf{H}_n)}\Big]^{1/2}\tau_n^{1/2}r\leq \delta$ for all $n \geq N_{\delta}^*$. By (\ref{eqD}),
$P(\|\widehat{\beta}_n - \beta_0\| \le \delta) \ge P(\Omega_n^*(r) \cap E_n^*(r)) \ge 1-\varepsilon$ for all $n \ge N_{\varepsilon,\delta}=\max(N_\varepsilon,N_\delta^*)$.
\end{proof}

\begin{remark}
{\rm
Since $\tau_n \leq m \widetilde{\lambda}_n$, Theorem \ref{existence-th} remains valid if we replace $\tau_n$ by $m \widetilde{\lambda}_n$ in conditions $(I_w^*)$ and $(L_w^*)$.}
\end{remark}

For the asymptotically normality of $\widehat{\beta}_n$, we consider the following condition:
\vspace{2mm}

\noindent $(CC)$ For any $r>0$ and $\delta >0$,
$$P(\sup_{\|\x\|=1}\sup_{\|\y\|=1}\sup_{\beta \in B_n^*(r)}|\x^T\textbf{H}_n^{-1/2}\mathcal{D}_n(\beta)\textbf{H}_n^{-1/2}\y-\x^T\y| \le \delta) \to 1.$$

\begin{lemma}
Condition $(CC)$ implies condition $(L_w^*)$.
\end{lemma}

\begin{proof}
We denote by $\Omega_n(\delta,r)$ the event in condition $(CC)$. Choose $\delta \in (0,1)$ arbitrary. In particular, on the event $\Omega_n(\delta,r)$, for any $\beta \in B_n^*(r)$ and for any $\x \in \mathbb{R}^p$ with $\|\x\|=1$,
$$|\x^T\H_n^{-1/2}\mathcal{D}_n(\beta)\H_n^{-1/2}\x-1| \le \delta,$$
and hence, $\x^T\H_n^{-1/2}\mathcal{D}_n(\beta) \H_n^{-1/2}\x\geq 1-\delta=:c_0$. Therefore \
\begin{equation}
\label{inclusion}
\Omega_n(\delta,r) \subset \Omega_n^*(r),
\end{equation}
where $\Omega_n^*(r)$ is the same event as in the proof of Theorem \ref{existence-th}.
Since $P(\Omega_n(\delta,r)) \to 1$, it follows that $P(\Omega_n^*(r)) \to 1$.
\end{proof}

Let $c_n=\lambda_{\rm max}(\textbf{M}_n^{-1}\textbf{H}_n)$. We consider the following boundedness condition:

\vspace{3mm}

$(B)$ There exists $c > 0$ such that $\tau_nc_n \le c$ for all $n$.

\begin{theorem}\label{theorem232}
Under conditions $(I_w^*)$, $(D_w^*)$, $(CC)$ and $(B)$,
$$\M_n^{-1/2}\g_n = \M_n^{-1/2}\H_n(\widehat{\beta}_n - \beta_0) + o_p(1).$$
\end{theorem}

\begin{proof} Let $\eta > 0$ and $\varepsilon > 0$ be arbitrary. We have to prove that there exists an integer $N_{\eta,\varepsilon} \ge 1$ such that for all $n \ge N_{\eta,\varepsilon}$,
$$P(\|\textbf{M}_n^{-1/2}\g_n - \textbf{M}_n^{-1/2}\textbf{H}_n(\widehat{\beta}_n - \beta_0)\| \le \eta) \ge 1-\varepsilon.$$

Let $c_0 \in (0,1)$ be a constant which will be specified later, $\delta = 1-c_0$ and $r=\frac{1}{c_0}\sqrt{\frac{2p}{\varepsilon}}$. Let $E_n^*(r)$, $\Omega_n^*(r)$ and $\widetilde{\Omega}_n(r)$ be the same events as in the proof of Theorem \ref{existence-th}. Let $\Omega_n(\delta,r)$ be the event in condition (CC). By (\ref{Omega-E}) and (\ref{inclusion}),
$$E_n^*(r) \cap \Omega_n^*(\delta,r) \subset E_n^*(r) \cap \Omega_n^*(r) \subset \widetilde{\Omega}_n(r).$$ Using Taylor's formula and the fact that $\g_n(\widehat{\beta}_n)=0$, we infer that there exists $\overline{\beta}_n \in B_n^*(r)$ such that
$$\g_n=\mathcal{D}_n(\overline{\beta}_n)(\widehat{\beta}_n - \beta_0).$$
We multiply this identity by the matrix $\textbf{M}_n^{-1/2}$. We denote $\textbf{U}_n(\beta)=\textbf{H}_n^{-1/2}\mathcal{D}_n(\beta)\textbf{H}_n^{-1/2}-\textbf{I}$. We obtain
\begin{align*}
\textbf{M}_n^{-1/2}\g_n&=
\textbf{M}_n^{-1/2}\mathcal{D}_n(\overline{\beta}_n)(\widehat{\beta}_n-\beta_0)\\
&=\textbf{M}_n^{-1/2}\textbf{H}_n^{1/2}\textbf{U}_n(\beta)\textbf{H}_n^{1/2}(\widehat{\beta}_n-\beta_0)+\textbf{M}_n^{-1/2}\textbf{H}_n(\widehat{\beta}_n-\beta_0).
\end{align*}

\noindent Note that
\begin{align*}
\|\textbf{M}_n^{-1/2}\textbf{H}_n^{1/2}\|^2&=\lambda_{\rm max}(\textbf{H}_n^{1/2}\textbf{M}_n^{-1}\textbf{H}_n^{1/2})=\frac{1}{\lambda_{\rm min}(\textbf{H}_n^{-1/2}\textbf{M}_n\textbf{H}_n^{-1/2})}\\&=\frac{1}{\lambda_{\rm min}(\textbf{H}_n^{-1}\textbf{M}_n)}=\lambda_{\rm max}(\textbf{M}_n^{-1}\textbf{H}_n)=c_n.
\end{align*}

\noindent On the event $E_n^*(r) \cap \Omega_n(\delta,r)$,  for all $\beta \in B_n^*(r)$, $\|\textbf{U}_n(\beta)\| \le c_1 \|\textbf{U}_n(\beta)\|_E \leq c_1p\delta$ (since all the elements of $\textbf{U}_n(\beta)$ are bounded in modulus by $\delta$), and so by condition $(B)$
\begin{align*}
\|\textbf{M}_n^{-1/2}\g_n - \textbf{M}_n^{-1/2}\textbf{H}_n(\widehat{\beta}_n-\beta)\| & \le (c_n\tau_n)^{1/2}c_1p\delta r \le c^{1/2}c_1p\delta r \\&= c^{1/2}c_1p \frac{1-c_0}{c_0} \sqrt{\frac{2p}{\varepsilon}} \le \eta,
\end{align*}

\noindent if we choose $c_0 \in (0,1)$ (depending on $\eta$ and $\varepsilon$) such that $\frac{1}{c_0} \le 1 + \frac{\eta}{c^{1/2}c_1p}\sqrt{\frac{\varepsilon}{2p}}$. This means that
$$E_n^*(r) \cap \Omega_n(\delta,r) \subset \{\|\textbf{M}_n^{-1/2}\g_n-\textbf{M}_n^{-1/2}\textbf{H}_n(\widehat{\beta}_n-\beta_0)\| \le \eta\}.$$
Similarly to (\ref{eqD}), it can be shown that there exists $N_{\eta,\varepsilon} \ge 1$ such that $$P(E_n^*(r) \cap \Omega_n(\delta,r)) \ge 1-\varepsilon$$ for all $n \ge N_{\eta,\varepsilon}$. The conclusion follows.
\end{proof}

We define $\widetilde{\textbf{Y}}_i=(A_i^*)^{-1/2}\varepsilon_i^*$ and $\gamma_n^{(D)}=\max\limits_{i \le n}\gamma_{n,i}^{(D)}$, where $$\gamma_{n,i}^{(D)}=\lambda_{\max}(\textbf{H}_n^{-1/2}\textbf{D}_i^T(\textbf{V}_i^*)^{-1}\textbf{D}_i\textbf{H}_n^{-1/2}).$$


We consider the following condition:
\vspace{3mm}

\noindent $(N_{\delta})$ There exist constants $\delta > 0$ and $K > 0$ such  that $\max_{i\leq n}\max_{j\leq m}E(\widetilde{Y}_{ij}^{2+2/\delta}|\X_i) \le K$, and there exists a constant $K_n^{(D)}>0$ such that  $\gamma_{n}^{(D)} \le K_n^{(D)}$ and $(c_n\widetilde{\lambda}_n)^{1+\delta}K_n^{(D)} \to 0$.

\vspace{3mm}

The following result is the analogue of Lemma 2 in \cite{XY} in our case.

\begin{theorem}\label{theorem233}
Under condition $(N_\delta)$, we have $$\M_n^{-1/2}\g_n \stackrel{d}{\to} N_p({\bf 0},\textbf{I}).$$
\end{theorem}
\begin{proof}
By the Cramer-Wold theorem, it suffices to show that for any $\lambda \in \mathbb{R}^p$ with $\|\lambda\|=1$,
\begin{equation}
\label{normal_conv1}
\lambda^T\M_n^{-1/2}\g_n \to N(0,1).
\end{equation}
Fix $\lambda \in \bR^p$ with $\|\lambda\|=1$. Then $\lambda^T\textbf{M}_n^{-1/2}\g_n=\displaystyle\sum_{i=1}^{n}Z_{n,i}$ where $Z_{n,i}=\lambda^T\textbf{M}_n^{-1/2}\textbf{D}_i^T(\V_i^*)^{-1}\varepsilon_i^*$. The variables $(Z_{n,i})_{i \le n}$ are independent. By Lemma \ref{lem1},
$$\mathrm{E}[Z_{n,i}]=\mathrm{E}[\lambda^T\textbf{M}_n^{-1/2}\textbf{D}_i^T(\V_i^*)^{-1}\mathrm{E}(\varepsilon_i^*|\textbf{X}_i)]=0.$$
Let $s_n^2=\mathrm{Var}(\lambda^T\textbf{M}_n^{-1/2}\g_n)$. By Lemma \ref{cov-gn-lemma},
\begin{align*}
s_n^2&=\displaystyle\sum_{i=1}^{n}\mathrm{E}(Z_{n,i}^2)
=\displaystyle\sum_{i=1}^{n}\mathrm{E}[\lambda^T\textbf{M}_n^{-1/2}\textbf{D}_i^T(\V_i^*)^{-1}\mathrm{E}[\varepsilon_i^*(\varepsilon_i^*)^T|\textbf{X}_i](\V_i^*)^{-1}\textbf{D}_i\textbf{M}_n^{-1/2}\lambda]
\\
&=\mathrm{E}[\lambda^T\textbf{M}_n^{-1/2}\Big(\displaystyle\sum_{i=1}^{n}
\textbf{D}_i^T(\V_i^*)^{-1}\Sigma_i^* (\V_i^*)^{-1}\textbf{D}_i\Big)\textbf{M}_n^{-1/2}\lambda]=\mathrm{E}[\lambda^T\textbf{M}_n^{-1/2}\textbf{M}_n^*\textbf{M}_n^{-1/2}\lambda]
\\&=\lambda^T\textbf{M}_n^{-1/2}\mathrm{E}[\textbf{M}_n^*]\textbf{M}_n^{-1/2}\lambda=1.
\end{align*}
\noindent By the Central Limit Theorem for triangular arrays (see e.g. Theorem 27.2 of \cite{billingsley95}), relation (\ref{normal_conv1}) will follow, once we prove that the following Lindeberg condition holds: for any $\e>0$,
\begin{equation}\label{lindeberg}
\lim\limits_{n \to \infty}\displaystyle\sum_{i=1}^n\mathrm{E}(Z_{n,i}^2{\rm I}\left\{\|Z_{n,i}\| \ge \varepsilon\right\})=0.
\end{equation}

It remains to prove (\ref{lindeberg}). Let $\e>0$ be arbitrary. Using the Cauchy-Schwartz inequality $\x^T\y \le \|\x\|\cdot\|\y\|$ for any $p$-dimensional vectors $\x$ and $\y$, we see that
\begin{align*}
Z_{n,i} & \le \|\lambda^T\textbf{M}_n^{-1/2}\D_i(\V_i^*)^{-1/2}\|^2 \cdot \|(\V_i^*)^{-1/2}\varepsilon_i^*\|^2 \\&=(\lambda^T\textbf{M}_n^{-1/2}\D_i^T(\V_i^*)^{-1}\D_i\textbf{M}_n^{-1/2}\lambda) \cdot ((\varepsilon_i^*)^T(\V_i^*)^{-1}\varepsilon_i^*)
\\&=\overline{\gamma}_{n,i}\cdot(\varepsilon_i^*)^T(\textbf{A}_i^*)^{-1/2}\textbf{R}_i^{-1}(\alpha)(\textbf{A}_i^*)^{-1/2}\varepsilon_i^*,
\end{align*}
\noindent where $\overline{\gamma}_{n,i}=\lambda\textbf{M}_n^{-1/2}\textbf{D}_i^T(\V_i^*)^{-1}\textbf{D}_i\textbf{M}_n^{-1/2}\lambda$. Using inequality \eqref{ineq}, it follows that
$$(\varepsilon_i^*)^T(\textbf{A}_i^*)^{-1/2}\textbf{R}_i^{-1}(\alpha)(\textbf{A}_i^*)^{-1/2}\varepsilon_i^* \le \lambda_{\max}(\textbf{R}_i^{-1}(\alpha))\|\widetilde{\textbf{Y}}_i\|^2 \le \widetilde{\lambda}_n\|\widetilde{\textbf{Y}}_i\|^2.$$
\noindent We obtain that
\begin{equation}\label{step1}
Z_{n,i}^2 \le \overline{\gamma}_{n,i}\widetilde{\lambda}_n\|\widetilde{\textbf{Y}}_i\|^2, \quad \mbox{for any} \ i\leq n.
\end{equation}

\noindent We also need another upper bound for $Z_{n,i}^2$, which is obtained as follows. By \eqref{ineq},
\begin{align*}
\overline{\gamma}_{n,i}& \leq
\lambda_{\max}(\textbf{H}_n^{-1/2}\D_i^T(\V_i^*)^{-1}\D_i\textbf{H}_n^{-1/2})\cdot\lambda^T\textbf{M}_n^{-1/2}\textbf{H}_n\textbf{M}_n^{-1/2}\lambda
\\& \le \gamma_{n,i}^{(D)}\lambda_{\max}(\textbf{M}_n^{-1/2}\textbf{H}_n\textbf{M}_n^{-1/2}).
\end{align*}
Recalling that $c_n=\lambda_{\max}(\textbf{M}_n^{-1}\textbf{H}_n)=\lambda_{\max}(\textbf{M}_n^{-1/2}\textbf{H}_n\textbf{M}_n^{-1/2})$, we obtain that $\overline{\gamma}_{n,i} \le \gamma_{n,i}^{(D)}c_n$, and hence
\begin{equation}\label{step2}
Z_{n,i}^2 \le c_n\widetilde{\lambda}_n\gamma_{n}^{(D)}\|\widetilde{\textbf{Y}}_i\|^2 , \quad \mbox{for any} \ i\leq n.
\end{equation}
\noindent Coming back to (\ref{lindeberg}), and using (\ref{step1}) and (\ref{step2}), we obtain:
\begin{align*}
\displaystyle\sum_{i=1}^n\mathrm{E}[Z_{n,i}^2{\rm I} & \{Z_{n,i} \ge \varepsilon\}]
 \le \widetilde{\lambda}_n\displaystyle\sum_{i=1}^n\mathrm{E}\Big[\overline{\gamma}_{n,i} \ \mathrm{E}\Big[\|\widetilde{\textbf{Y}}_i\|^2{\rm I}\left\{\|\widetilde{\textbf{Y}}_i\|^2 \ge \frac{\varepsilon^2}{c_n\widetilde{\lambda}_n\gamma_n^{(D)}}\right\} | \textbf{X}_i\Big]\Big]
\\& \le \widetilde{\lambda}_nm\displaystyle\sum_{i=1}^n\mathrm{E}\Big[\overline{\gamma}_{n,i}\mathrm{E}\Big[\frac{\|\widetilde{\textbf{Y}}_i\|^2}{m}\frac{(\|\widetilde{\textbf{Y}}_i\|^2/m)^{1/\delta}}{(\varepsilon^2/(mc_n\widetilde{\lambda}_n\gamma_n^{(D)}))^{1/\delta}}I\left\{\|\widetilde{\textbf{Y}}_i\|^2 \ge \frac{\varepsilon^2}{c_n\widetilde{\lambda}_n\gamma_n^{(D)}}\right\}|\textbf{X}_i\Big]\Big]
\\& \le \widetilde{\lambda}_nm\Big(\frac{mc_n\widetilde{\lambda}_nK_n^{(D)}}{\varepsilon^2}
\Big)^{1/\delta}\displaystyle\sum_{i=1}^n\mathrm{E}\Big[\overline{\gamma}_{n,i}
\mathrm{E}\Big[\Big(\frac{\|\widetilde{\textbf{Y}}_i\|^2}{m}\Big)^{1+1/(\delta)}|\textbf{X}_i\Big]\Big].
\end{align*}

\noindent Since the function $\phi(t)=t^{1+1/\delta}$ is convex,
$$\Big(\frac{1}{m}\|\widetilde{\textbf{Y}}_i\|^2\Big)^{1+1/\delta}
=\Big(\frac{1}{m}\displaystyle\sum_{j=1}^m\widetilde{Y}_{ij}^2\Big)^{1+1/\delta}
\le \frac{1}{m}\displaystyle\sum_{j=1}^m\widetilde{Y}_{ij}^{2+2/\delta},$$
\noindent and hence, by condition $(N_{\delta})$,
$$\mathrm{E}\Big[\Big(\frac{\|\widetilde{\textbf{Y}}_i\|^2}{m}\Big)^{1+1/\delta}|\textbf{X}_i\Big]
\le \frac{1}{m}\displaystyle\sum_{j=1}^m\mathrm{E}(\widetilde{Y}_{ij}^{2+2/\delta}|\textbf{X}_i) \le K.$$

\noindent Therefore
\begin{equation}\label{step3}
\displaystyle\sum_{i=1}^n\mathrm{E}[Z_{n,i}^2I\{Z_{n,i} \ge \varepsilon\}]
\le \widetilde{\lambda}_nm\Big(\frac{mc_n\widetilde{\lambda}_nK_n^{(D)}}{\varepsilon^2}\Big)^{1/\delta}K\mathrm{E}\Big[\displaystyle\sum_{i=1}^n\overline{\gamma}_{n,i}\Big].
\end{equation}

\noindent Note that, by the definition of $\overline{\gamma}_{n,i}$,
\begin{align*}
\displaystyle\sum_{i=1}^n\overline{\gamma}_{n,i} &=
\lambda^T\textbf{M}_n^{-1/2}\Big(\displaystyle\sum_{i=1}^n\textbf{D}_i^T(\textbf{V}_i^*)^{-1}\textbf{D}_i^T\Big)\textbf{M}_n^{-1/2}\lambda
=\lambda^T\textbf{M}_n^{-1/2}\textbf{H}_n^*\textbf{M}_n^{-1/2}\lambda.
\end{align*}

\noindent Taking expectation on both sides of the previous equality, we obtain:
\begin{align*}
\mathrm{E}\Big[\displaystyle\sum_{i=1}^n\overline{\gamma}_{n,i}\Big]
& =\lambda^T\textbf{M}_n^{-1/2}\mathrm{E}(\textbf{H}_n^*)\textbf{M}_n^{-1/2}\lambda
= \lambda^T\textbf{M}_n^{-1/2}\textbf{H}_n\textbf{M}_n^{-1/2}\lambda\leq c_n.
\end{align*}

\noindent Introducing this in (\ref{step3}), we obtain:
\begin{align*}
\displaystyle\sum_{i=1}^n\mathrm{E}[Z_{n,i}^2I{Z_{n,i} \ge \varepsilon}] & \le
\widetilde{\lambda}_nm\Big(\frac{mc_n\widetilde{\lambda}_nK_n^{(D)}}{\varepsilon^2}\Big)^{1/\delta}Kc_n = m\Big(\frac{m}{\varepsilon^2}\Big)^{1/\delta}(c_n\widetilde{\lambda}_n)^{1+1/\delta}(K_n^{(D)})^{1/\delta}.
\end{align*}
The last term converges to $0$ by condition $(N_{\delta})$. This finishes the proof of \eqref{lindeberg}.
\end{proof}

the following result is an immediate consequence of Theorems \ref{theorem232} and  \ref{theorem233}.

\begin{corollary}
\label{asy-normality}
Under conditions $(I_w^*)$, $(D_w^*)$, $(CC)$, $(B)$ and $(N_{\delta})$,
$$\M_n^{-1/2}\H_n(\widehat{\beta}_n-\beta) \to N_p({\bf 0},\textbf{I}).$$
\end{corollary}

\begin{remark}
\label{practice-remark}
{\rm In practice, we replace the matrices $\textbf{M}_n$ and $\textbf{H}_n$ by
\begin{align*}
& \widehat{\textbf{M}}_n = \displaystyle\sum_{i=1}^n\textbf{D}_i(\widehat{\beta}_n)^T\textbf{V}_i^*(\widehat{\beta}_n,\alpha)^{-1}\widehat{\Sigma}_i^*(\widehat{\beta}_n)\textbf{V}_i^*(\widehat{\beta}_n,\alpha)^{-1}\textbf{D}_i(\widehat{\beta}_n)
\\& \widehat{\textbf{H}}_n = \displaystyle\sum_{i=1}^n\textbf{D}_i(\widehat{\beta}_n)^T\textbf{V}_i^*(\widehat{\beta}_n,\alpha)^{-1}\textbf{D}_i(\widehat{\beta}_n),
\end{align*}

\noindent where $\widehat{\Sigma}_i^*(\beta)=(\textbf{Y}_i^*-\mu_i(\beta))(\textbf{Y}_i^*-\mu_i(\beta))^T$. Note that the weighted response $Y_{ij}^*$ depends on the missingness probability $\pi_{ij}$ which is unknown (see definition \eqref{def-Ystar} of $Y_{ij}^*$). Moreover, the matrix $\textbf{V}_i^*(\widehat{\beta}_n,\alpha)$ depends on $\textbf{A}_i^*(\widehat{\beta}_n)$ (see (\ref{def_Vi_star})), which also depends on the probabilities ${(\pi_{ij})}_{1 \le j \le m}$ (see (\ref{sigma_star_ijj}) for the definition of components $\sigma_{i,jj}^*(\beta),j=1,\ldots,m$ of the diagonal matrix $\As_i(\beta)$). To avoid this problem, we may use a logistic regression model to ``posit the missing data process'', as suggested on page 155 of \cite{yi-ma-carroll12}. This consists in fitting a logistic regression model to the complete data consisting of $(\textbf{I}_i,\textbf{X}_i)$ for $i=1,\dots,n$, with a new regression parameter $\gamma$. As in Example \ref{example3} (applied to the case when $Y_{ij}$ is replaced by $I_{ij}$), we assume that $I_{ij}$ is a Bernoulli random variable with mean
\begin{equation}
\label{def-pi-ij}
\pi_{ij}=\pi_{ij}(\gamma)=\frac{\exp(\X_{ij}^T\gamma)}{1+\exp(\X_{ij}^T\gamma)}.
\end{equation}
To estimate $\gamma$, we solve the classical GEE with working independence matrices $\textbf{R}_i(\alpha)=\textbf{I}$ for all $i=1,\dots,n$: (see equation (\ref{twooneoneprime}))
\begin{equation}\label{eq50}
\displaystyle\sum_{i=1}^n\displaystyle\sum_{j=1}^m \X_{ij}\Big(I_{ij}-\frac{\exp(\X_{ij}^T\gamma)}{1+\exp(\X_{ij}^T\gamma)}\Big)=0.
\end{equation}
Let $\widehat{\gamma}$ be the solution of the equation (\ref{eq50}). Then, in the calculation of $Y_{ij}^*$ and $\sigma_{i,jj}^*(\widehat{\beta}_n)$, we replace $\pi_{ij}$ by $\widehat{\pi}_{ij}=\pi_{ij}(\widehat{\gamma})$. }
\end{remark}

\section{Verification of condition (CC)}
\label{section-CC}

In this section, we give some sufficient conditions which ensure that condition (CC) holds.  Proceeding as in Remark 1 of \cite{XY} (see also Appendix A of \cite{XY}), we write the derivative of $\g_n(\beta)$ as the sum of three terms:
$$-\cD_n(\beta)=\frac{\partial }{\partial \beta^T}\g_n(\beta)=-\Hs_n(\beta)+\B_n(\beta)+\cE_n(\beta),$$
where $\Hs_n(\beta)$ is given by \eqref{def-Hstar}, $\B_n(\beta)=\B_n^{(1)}(\beta)+\B_n^{(2)}(\beta)$, $\cE_n(\beta)=\cE_n^{(1)}(\beta)+\cE_n^{(2)}(\beta)$ and
\begin{align*}
\B_n^{(1)}(\beta)&=\sum_{i=1}^{n}\X_i^T {\rm diag}\big[\Ri^{-1}\As_i(\beta)^{-1/2}(\mu_i-\mu_i(\beta))\big] \G_i^{(1)}(\beta)\X_i \\
\B_n^{(2)}(\beta)&=\sum_{i=1}^{n}\X_i^T \F_i(\beta) \Ri^{-1} {\rm diag}\big[\mu_i-\mu_i(\beta)\big] \G_i^{(2)}(\beta)\X_i\\
\cE_n^{(1)}(\beta)&=\sum_{i=1}^{n}\X_i^T {\rm diag}\big[\Ri^{-1}\As_i(\beta)^{-1/2}\es_i \big] \G_i^{(1)}(\beta)\X_i \\
\cE_n^{(2)}(\beta)&=\sum_{i=1}^{n}\X_i^T \F_i(\beta) \Ri^{-1} {\rm diag}\big[\es_i\big] \G_i^{(2)}(\beta)\X_i.
\end{align*}

Here $\G_i^{(k)}(\beta)={\rm diag}\big(g_{i1}^{(k)}(\beta),\ldots,g_{im}^{(k)}(\beta)\big)$ for $k=1,2$, where
$$\frac{\partial}{\partial \beta^T} f_{ij}(\beta)=g_{ij}^{(1)}(\beta)\X_{ij}^T \quad \mbox{and} \quad \frac{\partial}{\partial \beta^T} [\sigma_{i,jj}^{*}(\beta)]^{-1/2}=g_{ij}^{(2)}(\beta)\X_{ij}^T,$$
with functions $g_{ij}^{(1)}(\beta)$ and $g_{ij}^{(2)}(\beta)$ given by:
\begin{align}
\label{def-gij-1}
g_{ij}^{(1)}(\beta)&=\frac{\mu''(\X_{ij}^T \beta)}{[\sigma_{i,jj}^{*}(\beta)]^{1/2}}-\frac{2 \left(\frac{1}{\pi_{ij}}-1\right)\mu(\X_{ij}^T \beta) (\mu'(\X_{ij}^{T}\beta))^2+\frac{1}{\pi_{ij}} \mu'(\X_{ij}^T \beta)\mu''(\X_{ij}^T \beta)}{2[\sigma_{i,jj}^{*}(\beta)]^{3/2}}\\
\label{def-gij-2}
g_{ij}^{(2)}(\beta)&=-\frac{2 \left(\frac{1}{\pi_{ij}}-1\right)\mu(\X_{ij}^T \beta) \mu'(\X_{ij}^{T}\beta)+\frac{1}{\pi_{ij}} \mu''(\X_{ij}^T \beta)}{2[\sigma_{i,jj}^{*}(\beta)]^{3/2}}
\end{align}

We treat separately the three terms. For this, we introduce the same constants and smoothness assumption as on pages 330-331 of \cite{XY}:

$$\gamma_n^{(0)}=\max_{i\leq n}\max_{j \leq m} (\X_{ij}^T \H_n^{-1}\X_{ij}), \qquad \gamma_n^{*}=\tau_n \gamma_n^{(0)}, \qquad \pi_n=\frac{\max_{1 \leq i \leq n} \lM(\Ri^{-1})}{\min_{1\leq i\leq n}\lm(\Ri^{-1})}.$$

\vspace{3mm}

{\bf Assumption (AH).} $k_n^{(i)}=O_p(1)$ for $i=1,2,3$
where
$$k_n^{(0)}=\sup_{\beta \in B_n^*(r)}\max_{i,j}\frac{\mu'(\X_{ij}^{T}\beta)}{\mu(\X_{ij}^{T}\beta)} \quad \qquad k_n^{(1)}=\sup_{\beta \in B_n^*(r)}\max_{i,j}\frac{\mu''(\X_{ij}^{T}\beta)}{\mu'(\X_{ij}^{T}\beta)}$$

We impose the following assumption on the missingness probabilities:
\vspace{3mm}

{\bf Assumption (M).} $\rho_n=O_p(1)$, where
$$\rho_n=\max_{i\leq n}\max_{j\leq m}\frac{1}{\pi_{ij}}.$$

Assumption M says that for any $\e>0$, there exists a constant $C_{\e}>0$ and an integer $N_{\e}\geq 1$ such that for any $n\geq N_{\e}$, with probability greater than $1-\e$, $\pi_{ij} \geq C_{\e}$ for all $i\leq n$ and $j\leq m$. Intuitively speaking, this means that the missingness probabilities $\pi_{ij}$ are bounded away from $0$. Note that the case when all probabilities $\pi_{ij}$ are equal to $0$ corresponds to the case when all the data is missing.

The following three lemmas are the counterparts of Lemmas A.1.(ii), A.2.(ii) and A.3.(ii) of \cite{XY}, when the covariates are random and the responses are missing at random.

\begin{lemma}
\label{lemmaA1}
Suppose Assumptions (AH) and (M) hold. If $\pi_n \gamma_n^{*} \stackrel{P}{\to} 0$ then
$$\sup_{\|\x\|=1} \sup_{\|\y\|=1} \sup_{\beta \in B_n^{*}(r)}|\x^{T} \H_n^{-1/2}\Hs_n(\beta)\H_n^{-1/2}\y-\x^{T}\y| \stackrel{P}{\to}0.$$
\end{lemma}

\begin{lemma}
\label{lemmaA2}
Suppose Assumptions (AH) and (M) hold. If $\pi_n^2 \gamma_n^{*} \stackrel{P}{\to} 0$ then
$$\sup_{\|\x\|=1} \sup_{\|\y\|=1} \sup_{\beta \in B_n^{*}(r)}|\x^{T} \H_n^{-1/2}\B_n(\beta)\H_n^{-1/2}\y| \stackrel{P}{\to}0.$$
\end{lemma}

\begin{lemma}
\label{lemmaA3}
Suppose Assumptions (AH) and (M) hold. If $\gamma_{n}^* \stackrel{P}{\to}0$, $n\pi_n^2 \gamma_n^{(0)}=O_p(1)$ and $n\pi_n^2 \gamma_n^{(0)}E(\tau_n)\stackrel{P}{\to} 0$
then
$$\sup_{\|\x\|=1} \sup_{\|\y\|=1} \sup_{\beta \in B_n^{*}(r)}|\x^{T} \H_n^{-1/2}\cE_n(\beta)\H_n^{-1/2}\y| \stackrel{P}{\to}0.$$
\end{lemma}

\noindent {\bf Proof of Lemma \ref{lemmaA1}:} Writing $\F_i(\beta)=\F_i+(\F_i(\beta)-\F_i)$ in definition \eqref{def-Hstar}, we obtain that
$$\x^{T} \H_n^{-1/2}\Hs_n(\beta)\H_n^{-1/2}\y-\x^T \y = T_0(\x,\y)+\sum_{i=1}^{3}T_i(\beta,\x,\y),$$
where $T_0(\x,\y)=\x^{T} \H_n^{-1/2}\Hs_n\H_n^{-1/2}\y-\x^T \y$ and
\begin{align}
\label{def-T1}
T_1(\beta,\x,\y)&=\sum_{i=1}^{n}\x^T \H_n^{-1/2} \X_i^T (\F_i(\beta)-\F_i)\Ri^{-1}
(\F_i(\beta)-\F_i) \X_i \H_n^{-1/2}\y\\
\nonumber
T_2(\beta,\x,\y)&=\sum_{i=1}^{n}\x^T \H_n^{-1/2} \X_i^T (\F_i(\beta)-\F_i)\Ri^{-1}
\F_i\X_i \H_n^{-1/2}\y\\
\nonumber
T_3(\beta,\x,\y)&=\sum_{i=1}^{n}\x^T \H_n^{-1/2} \X_i^T \F_i\Ri^{-1}
(\F_i(\beta)-\F_i)\X_i \H_n^{-1/2}\y.
\end{align}

To treat $T_0(\x,\y)$, note that $\Hs_n=\sum_{i=1}^{n}\U_{i}$, where $\U_i=\D_i^T (\Vs_i)^{-1}\D_i, i=1,\ldots,n$ are i.i.d. random matrices. By the strong law of large numbers, $\frac{1}{n}\Hs_n \to E(\U_1)$ a.s. (component-wise), and hence
$\|\frac{1}{n}\Hs_n -E(\U_1)\|\to 0$ a.s. Since $\H_n=nE(\U_1)$, we obtain:
\begin{equation}
\label{Hn-conv}
\|\H_n^{-1/2}\Hs_n \H_n^{-1/2}-\I\| \to 0 \quad \mbox{a.s.}
\end{equation}
Therefore, $\sup_{\x,\y}|T_0(\x,\y)| \to 0$ a.s. Using inequality \eqref{ineq}, we have:
\begin{align}
\nonumber
\sup_{\x,\y,\beta}|T_1(\beta,\x,\y)| &\leq \pi_n \sup_{\beta}\max_{i \leq n}\lM^2(\F_i^{-1}\F_i(\beta)-\I)\cdot \sup_{\x,\y}|\x^T \H_n^{-1/2}\Hs_n \H_n^{-1/2}\y|\\
\label{bound-T1}
&=\pi_n O_{p}(\gamma_n^{*})O_{p}(1)=\pi_n \gamma_n^{*}O_{p}(1)=o_p(1),
\end{align}
where the first equality above is due to  Lemma \ref{lemmaF} (Appendix A) and relation \eqref{Hn-conv}.

To treat $T_2(\beta,\x,\y)$, we use Cauchy-Schwatz inequality: for any $p$-dimensional vectors $(\a_i)_{i=1,\ldots,n}$ and $(\b_i)_{i=1,\ldots, n}$,
\begin{equation}
\label{CS}
\left|\sum_{i=1}^{n}\a_i^T \b_i \right| \leq \Big(\sum_{i=1}^{n}\a_i^T \a_i\Big)^{1/2} \Big(\sum_{i=1}^{n}\b_i^T \b_i\Big)^{1/2}.
\end{equation}
Letting $\a_i=\x^T \H_n^{-1/2} \X_i^T \F_i\Ri^{-1/2}$ and $\b_i^T=\Ri^{1/2}
(\F_i^{-1}\F_i(\beta)-\I)\Ri^{-1}\F_i\X_i \H_n^{-1/2}\y$, we obtain
$|T_2(\beta,\x,\y)| \leq T_2'(\x)^{1/2} T_2''(\beta,\y)^{1/2}$,
with $T_2'(\x)=\x \H_n^{-1/2}\Hs_n \H_n^{-1/2}\x$ and
\begin{align*}
T_2''(\beta,\y)& \leq \pi_n \max_{i \leq n}\lM^2(\F_i^{-1} \F_i(\beta)-\I)\y \H_n^{-1/2}\Hs_n \H_n^{-1/2}\y.
\end{align*}
Arguing as above, we get $\sup_{\beta,\x,\y}|T_2(\beta,\x,\y)|=o_{p}(1)$. The term $T_3(\beta,\x,\y)$ is similar. $\Box$

\vspace{3mm}

\noindent {\bf Proof of Lemma \ref{lemmaA2}:} We begin by treating $\B_n^{(1)}(\beta)$. Note that for any $p\times p$ diagonal matrix $\Delta$ and for any $p$-dimensional vectors $\v$ and $\w$,
\begin{equation}
\label{diagonal}
{\rm diag}(\v)\Delta \w=\Delta {\rm diag}(\w)\v.
\end{equation}
We use this with $\v=\Ri^{-1}\As_i(\beta)^{-1/2}(\mu_i-\mu_i(\beta))$,
$\Delta=\G_i^{(1)}(\beta)$ and $\w=\X_i \H_n^{-1/2}\y$. We obtain that $\x^{T} \H_n^{-1/2}\B_n^{(1)}(\beta)\H_n^{-1/2}\y$ is equal to
$$\sum_{i=1}^{n}\x^T \H_n^{-1/2}\X_i^T \G_i^{(1)}(\beta){\rm diag}(\X_i \H_n^{-1/2}\y)\Ri^{-1}\As_i(\beta)^{-1/2}(\mu_i-\mu_i(\beta)).$$
Using Cauchy-Schwarz inequality \eqref{CS}, it follows that
\begin{equation}
\label{bound-B1}
|\x^{T} \H_n^{-1/2}\B_n^{(1)}(\beta)\H_n^{-1/2}\y| \leq S_1(\beta,\x,\y)^{1/2} S_2(\beta)^{1/2},
\end{equation}
where
\begin{align}
\label{def-S1}
S_1(\beta,\x,\y)&=\sum_{i=1}^{n}\x^T \H_n^{-1/2}\X_i^T \G_i^{(1)}(\beta){\rm diag}(\X_i \H_n^{-1/2}\y)\Ri^{-1}{\rm diag}(\X_i \H_n^{-1/2}\y)\\
\nonumber
& \quad \quad \quad \G_i^{(1)}(\beta)\X_i \H_n^{-1/2}\x\\
\label{def-S2}
S_2(\beta)&=\sum_{i=1}^{n} (\mu_i-\mu_i(\beta))^T \As_i(\beta)^{-1/2} \Ri^{-1}\As_i(\beta)^{-1/2}(\mu_i-\mu_i(\beta)).
\end{align}
Using \eqref{ineq} and the fact that
\begin{equation}
\label{diag-Xi}
\lM^2({\rm diag}(\X_i \H_n^{-1/2}\y))\leq \gamma_n^{(0)},
\end{equation}
it follows that
\begin{align}
\label{bound-S1}
& S_1(\beta,\x,\y) \leq \tilde{\lambda}_n \gamma_n^{(0)} \sum_{i=1}^{n}\x^T \H_n^{-1/2}\X_i^T [\G_i^{(1)}(\beta)]^2 \X_i \H_n^{-1/2}\x \\
\nonumber
& \leq \tilde{\lambda}_n \gamma_n^{(0)}  \max_{i\leq n}
\lM \big(\Ri^{1/2} \F_i^{-1}[\G_i^{(1)}(\beta)]^2 \F_i^{-1}   \Ri^{1/2} \big) \cdot \x^T \H_n^{-1/2}\Hs_n \H_n^{-1/2}\x \\
\nonumber
&\leq \pi_n \gamma_n^{(0)} \max_{i \leq n}\lM^2\big(\F_i^{-1} \G_i^{(1)}(\beta) \big)
\cdot \x\H_n^{-1/2}\Hs_n \H_n^{-1/2}\x.
\end{align}
By relations \eqref{fij-Op} and \eqref{gij1-fij} (given in Appendix A),  $\max_{i \leq n}\lM^2\big(\F_i^{-1}\G_i^{(1)}(\beta)\big)=O_p(1)$. By Lemma \ref{lemmaA1},
$\sup_{\|\x\|=1}\big(\x\H_n^{-1/2}\Hs_n \H_n^{-1/2}\x\big)=O_p(1)$.
From this, we infer that
\begin{equation}
\label{S1-beta-xy}
\sup_{\|\x\|=1} \sup_{\|\y\|=1} \sup_{\beta \in B_n^*(r)}S_1(\beta,\x,\y) \leq \pi_n \gamma_n^{(0)}O_p(1).
\end{equation}

We now treat $S_2(\beta)$. By Taylor's formula, for any $\beta \in B_n^*(r)$, there exists $\overline{\beta}_{ij} \in B_n^*(r)$ such that $\mu_{ij}(\beta)-\mu_{ij}(\beta_0)=\mu'(\X_{ij}^T\overline{\beta}_{ij})\X_{ij}^T (\beta-\beta_0)$. Then
$\mu_i(\beta)-\mu_i=\Ab_i\X_i(\beta-\beta_0)$, where
$\Ab_i$ is the diagonal matrix with entries $\mu'(\X_{ij}^T \overline{\beta}_{ij}),j=1,\ldots,m$. Note that $\Ab_i \As_i(\beta)^{-1/2} =\As_i(\beta)^{-1/2} \Ab_i $ since $\Ab_i$ and $\As_i(\beta)^{-1/2}$ are diagonal matrices. Using inequality \eqref{ineq}, we get:
\begin{align*}
S_2(\beta)&=\sum_{i=1}^{n} (\beta-\beta_0)^T\X_i^T \As_i(\beta)^{-1/2} \Ab_i \Ri^{-1}\Ab_i \As_i(\beta)^{-1/2}\X_i(\beta-\beta_0)\\
&=\sum_{i=1}^{n} (\beta-\beta_0)^T\X_i^T \F_i(\beta)\A_i(\beta)^{-1}\Ab_i \Ri^{-1}\Ab_i \A_i(\beta)^{-1}\F_i(\beta)\X_i(\beta-\beta_0)\\
&\leq \widetilde{\lambda}_n \max_{i\leq n}\lM^2(\Ab_i \A_i(\beta)^{-1}) \sum_{i=1}^{n}(\beta-\beta_0)^T \X_i^T [\F_i(\beta)]^2\X_i (\beta-\beta_0)\\
&\leq \pi_n \max_{i\leq n}\lM^2(\Ab_i \A_i(\beta)^{-1}) \cdot (\beta-\beta_0)^T \H_n^*(\beta)(\beta-\beta_0)\\
&\leq  \pi_n \max_{i\leq n}\lM^2(\Ab_i \A_i(\beta)^{-1}) \cdot \|\H_n^{-1/2}\Hs_n(\beta) \H_n^{-1/2}\| \cdot \|\H_n^{-1/2}(\beta-\beta_0)\|^2\\
&\leq  \pi_n \max_{i\leq n}\lM^2(\Ab_i \A_i(\beta)^{-1}) \cdot \|\H_n^{-1/2}\Hs_n(\beta) \H_n^{-1/2}\| \cdot \tau_n r^2
\end{align*}
By Lemma \ref{lemmaA1}, $\sup_{\beta \in B_n^*(r)}\|\H_n^{-1/2}\Hs_n(\beta) \H_n^{-1/2}\|=O_p(1)$. Note that $\Ab_i \A_i(\beta)^{-1}$ is a diagonal matrix with entries $\mu'(\X_{ij}^T \beta_{ij})/\mu'(\X_{ij}^T \beta),j=1,\ldots,m$. By relation \eqref{fij-Op} (given in Appendix A), it follows that
$\sup_{\beta \in B_{n}^{*}(r)}\max_{i\leq n}\lM^2(\Ab_i \A_i(\beta)^{-1})=O_p(1)$. Hence,
\begin{equation}
\label{S2-beta}
\sup_{\beta \in B_n^*(r)}S_2(\beta) \leq \pi_n \tau_n O_p(1).
\end{equation}

Using relations \eqref{bound-B1}, \eqref{S1-beta-xy} and \eqref{S2-beta}, we infer that
$$\sup_{\|x\|=1}\sup_{\|y\|=1}\sup_{\beta \in B_n^{*}(r)}|\x^{T} \H_n^{-1/2}\B_n^{(1)}(\beta)\H_n^{-1/2}\y| \leq \pi_n (\gamma_n^{*})^{1/2} O_p(1)=o_p(1).$$

We continue with the treatment of $\B_n^{(2)}(\beta)$. Using relation \eqref{diagonal}, we see that
$$\x^{T} \H_n^{-1/2}\B_n^{(2)}(\beta)\H_n^{-1/2}\y=\sum_{i=1}\x^T \H_n^{-1/2}\X_i^T\F_i(\beta)\Ri^{-1}{\rm diag}(\X_i \H_n^{-1/2}\y)\G_i^{(2)}(\beta)(\mu_i-\mu_i(\beta)).$$
We use Cauchy-Schwarz inequality \eqref{CS} with $\a_i^T=\x^T \H_n^{-1/2}\X_i^T\F_i(\beta)\Ri^{-1}{\rm diag}(\X_i \H_n^{-1/2}\y)\linebreak \G_i^{(2)}(\beta)\As_i(\beta)^{1/2}\Ri^{1/2}$ and $\b_i=\Ri^{-1/2}\As_i(\beta)^{-1/2}(\mu_i-\mu_i(\beta))$. We obtain:
\begin{equation}
\label{bound-B2}
|\x^{T} \H_n^{-1/2}\B_n^{(2)}(\beta)\H_n^{-1/2}\y| \leq S_3(\beta,\x,\y)^{1/2}S_2(\beta)^{1/2},
\end{equation}
where $S_2(\beta)$ is given by \eqref{def-S2} and
\begin{align}
\label{def-S3}
S_3(\beta,\x,\y)&=\sum_{i=1}^{n}
\x^T \H_n^{-1/2}\X_i^T\F_i(\beta)\Ri^{-1}{\rm diag}(\X_i \H_n^{-1/2}\y) \G_i^{(2)}(\beta)\As_i(\beta)^{1/2}\Ri\\
\nonumber
& \quad \quad \quad  \As_i(\beta)^{1/2}\G_i^{(2)}(\beta){\rm diag}(\X_i \H_n^{-1/2}\y) \Ri^{-1}\F_i(\beta)\X_i \H_n^{-1/2}\x.
\end{align}
Using inequalities \eqref{ineq} and \eqref{diag-Xi}, we obtain that:$$
S_3(\beta,\x,\y)\leq \pi_n \gamma_n^{(0)}\max_{i \leq n}\lM^2(\As_i(\beta)^{1/2}\G_i^{(2)}(\beta))\cdot \x^T \H_n^{-1/2}\Hs_n(\beta)\H_n^{-1/2}\x.$$
Note that $\As_i(\beta)^{1/2}\G_i^{(2)}(\beta)$ is a diagonal matrix with elements $\sqrt{\sigma_{i,jj}^*(\beta)}g_{ij}^{(2)}(\beta),j=1,\ldots,m$. By Lemma \ref{lemmaD} (Appendix A), $\sup_{\beta \in B_n^*(r)}\max_{i \leq n}\lM^2(\As_i(\beta)^{1/2}\G_i^{(2)}(\beta))=O_p(1)$. Using Lemma \ref{lemmaA1}, we obtain:
\begin{equation}
\label{S3-beta-xy}
\sup_{\|x\|=1}\sup_{\|y\|=1}\sup_{\beta \in B_n^*(r)}S_3(\beta,\x,\y) \leq \pi_n \gamma_n^{(0)}O_p(1).
\end{equation}
Using relations \eqref{bound-B2}, \eqref{S3-beta-xy} and \eqref{S2-beta}, we infer that:
$$\sup_{\|x\|=1}\sup_{\|y\|=1}\sup_{\beta \in B_n^{*}(r)}|\x^{T} \H_n^{-1/2}\B_n^{(2)}(\beta)\H_n^{-1/2}\y| \leq \pi_n (\gamma_n^{*})^{1/2} O_p(1)=o_p(1).$$
$\Box$

\vspace{3mm}

\noindent {\bf Proof of Lemma \ref{lemmaA3}:} We first treat the term $\cE_n^{(1)}(\beta)$.
Using relation \eqref{diagonal}, we see that
\begin{align*}
\x^T \H_n^{-1/2} \cE_n^{(1)}(\beta)\H_n^{-1/2}\y&=\sum_{i=1}^{n}\x^T \H_n^{-1/2}\X_i^T \G_i^{(1)}(\beta){\rm diag}(\X_i \H_n^{-1/2}\y)\Ri^{-1} \As_i(\beta)^{-1/2}\e_i^*\\
&=U_1(\x,\y)+U_3(\beta,\x,\y)+U_5(\beta,\x,\y),
\end{align*}
where
\begin{align*}
U_1(\x,\y)&=\sum_{i=1}^{n}\x^T \H_n^{-1/2}\X_i^T \G_i^{(1)}{\rm diag}(\X_i \H_n^{-1/2}\y)\Ri^{-1} (\As_i)^{-1/2}\e_i^*\\
U_3(\beta,\x,\y)&=\sum_{i=1}^{n}\x^T \H_n^{-1/2}\X_i^T \G_i^{(1)}(\beta){\rm diag}(\X_i \H_n^{-1/2}\y)\Ri^{-1} \big(\As_i(\beta)^{-1/2}-(\As_i)^{-1/2}\big)\e_i^*\\
U_5(\beta,\x,\y)&=\sum_{i=1}^{n}\x^T \H_n^{-1/2}\X_i^T \big(\G_i^{(1)}(\beta)-\G_i^{(1)}\big) {\rm diag}(\X_i \H_n^{-1/2}\y)\Ri^{-1} (\As_i)^{-1/2}\e_i^*.
\end{align*}

We first treat $U_1(\x,\y)$. By the Cauchy-Schwarz inequality \eqref{CS},
\begin{equation}
\label{bound-U1}
|U_1(\x,\y)|\leq S_1(\beta_0,\x,\y)^{1/2}U^{1/2},
\end{equation}
where $S_1(\beta,\x,\y)$ is given by \eqref{def-S1} and
$U=\sum_{i=1}^{n}W_i$, with $W_i= (\e_i^*)^T (\As_i)^{-1/2} \Ri^{-1} \linebreak (\As_i)^{-1/2}\e_i^*$. Using the fact that $\x^T \x={\rm tr}(\x \x^T)$ for any $p$-dimensional vector $\x$, we obtain:
\begin{align*}
E(W_i)&= E[E[{\rm tr} \{\Ri^{-1/2} (\As_i)^{-1/2} \e_i^* (\e_i^*)^T (\As_i)^{-1/2} \Ri^{-1/2}\}|\X_i]\\
&=E[{\rm tr} \{\Ri^{-1/2} (\As_i)^{-1/2} E[\e_i^* (\e_i^*)^T|\X_i] (\As_i)^{-1/2} \Ri^{-1/2} \}]\\
&= E[{\rm tr} \{ \Ri^{-1/2} \Rs_i \Ri^{-1/2}\} ]\leq mE[\lM(\Ri^{-1/2} \Rs_i \Ri^{-1/2})]\leq mE(\tau_n),
\end{align*}
for any $i=1,\ldots,n$, using \eqref{sigma-R-star} for the last equality. Hence,
$\sum_{i=1}^{n}E(W_i) \leq mn E(\tau_n)$. Since $\{(\Y_i,\X_i,\I_i)\}_{i=1,\ldots,n}$ are i.i.d., $(W_i)_{i=1,\ldots,n}$ are independent. Therefore, by Chebyshev's weak law of large numbers,
$\sum_{i=1}^{n}(W_i-E(W_i))=o_p(n)$. Hence,
\begin{equation}
\label{bound-U}
U \leq o_p(n)+mn E(\tau_n).
\end{equation}

Using \eqref{bound-U1}, \eqref{S1-beta-xy}, \eqref{bound-U} and the hypotheses of the lemma, it follows that
$$\sup_{\|\x\|=1} \sup_{\|\y\|=1} |U_1(\x,\y)| \leq \{\big(\pi_n \gamma_n^{(0)}o_p(n)\big)^{1/2}+\big(\pi_n \gamma_n^{(0)}n E(\tau_n)\big)^{1/2}\} O_p(1)=o_p(1).$$

Next, we treat $U_3(\beta,\x,\y)$. By the Cauchy-Schwartz inequality \eqref{CS}, it follows that
\begin{equation}
\label{bound-U3}
|U_3(\beta,\x,\y)|\leq U_3'(\beta,\x,\y)^{1/2}U^{1/2},
\end{equation}
where $U$ is the same as above and
\begin{align*}
U_3'(\beta,\x,\y)&=\sum_{i=1}^{n}\x^T \H_n^{-1/2}\G_i^{(1)}(\beta){\rm diag}(\X_i \H_n^{-1/2}\y)\Ri^{-1}\big(\As_i(\beta)^{-1/2}(\As_i)^{1/2}-\I\big) \Ri\\
& \big(\As_i(\beta)^{-1/2}(\As_i)^{1/2}-\I\big) \Ri^{-1} {\rm diag}(\X_i \H_n^{-1/2}\y) \G_i^{(1)}(\beta) \H_n^{-1/2}\x.
\end{align*}
Using inequalities \eqref{ineq} and \eqref{diag-Xi}, we see that
\begin{align*}
U_3'(\beta,\x,\y)&\leq \pi_n \widetilde{\lambda}_n \gamma_n^{(0)}\max_{i\leq n}\lM^2(\As_i(\beta)^{-1/2}(\As_i)^{1/2}-\I)  \sum_{i=1}^{n}\x^T \H_n^{-1/2}\X_i^T [\G_i^{(1)}(\beta)]^2 \X_i \H_n^{-1/2}\x.
\end{align*}
Proceeding as in \eqref{bound-S1} and using Lemma \ref{lemmaC} (Appendix A), we get:
\begin{equation}
\label{bound-U3'}
\sup_{\|x\|=1}\sup_{\|y\|=1}\sup_{\beta \in B_n^*(r)}U_3'(\beta,\x,\y) \leq \pi_n^2 \gamma_n^{(0)}O_{p}(1).
\end{equation}
Using \eqref{bound-U3}, \eqref{bound-U3'} and \eqref{bound-U}, we obtain by the hypotheses of the lemma that
$$\sup_{\|x\|=1}\sup_{\|y\|=1}\sup_{\beta \in B_n^*(r)}|U_3(\beta,\x,\y)|\leq  \{\big(\pi_n^2 \gamma_n^{(0)}o_p(n)\big)^{1/2}+\big(\pi_n^2 \gamma_n^{(0)}n E(\tau_n)\big)^{1/2}\} O_p(1)=o_p(1).$$

We now treat $U_5(\beta,\x,\y)$. By the Cauchy-Schwartz inequality \eqref{CS}, it follows that
\begin{equation}
\label{bound-U5}
|U_5(\beta,\x,\y)|\leq U_5'(\beta,\x,\y)^{1/2}U^{1/2},
\end{equation}
where $U$ is the same as above and
\begin{align*}
U_5'(\beta,\x,\y)&=\sum_{i=1}^{n}\x^T \H_n^{-1/2}\big(\G_i^{(1)}(\beta)-\G_i^{(1)}\big){\rm diag}(\X_i \H_n^{-1/2}\y)\Ri^{-1}{\rm diag}(\X_i \H_n^{-1/2}\y) \\
&\quad \quad \quad \big(\G_i^{(1)}(\beta)-\G_i^{(1)}\big)\X_i \H_n^{-1/2}\x.
\end{align*}

Using inequalities \eqref{ineq} and \eqref{diag-Xi}, it follows that
\begin{align*}
U_5'(\beta,\x,\y)& \leq \pi_n \gamma_n^{(0)}\max_{i \leq n}\lM^2(\F_i^{-1} (\G_i^{(1)}(\beta)-\G_i^{(1)})) \cdot \x^T \H_n^{-1/2}\Hs_n \H_n^{-1/2}\x.
\end{align*}
The matrix $\F_i^{-1} (\G_i^{(1)}(\beta)-\G_i^{(1)})$ has $j$-th element given by
$$\frac{g_{ij}^{(1)}(\beta)-g_{ij}^{(1)}(\beta_0)}{f_{ij}(\beta_0)}=
\frac{g_{ij}^{(1)}(\beta)}{f_{ij}(\beta)}\cdot \frac{f_{ij}(\beta)}{f_{ij}(\beta_0)}-\frac{g_{ij}^{(1)}(\beta_0)}{f_{ij}(\beta_0)}.$$
By relation \eqref{gij1-fij} (Appendix A), $\max_{i \leq n}\lM^2(\F_i^{-1} (\G_i^{(1)}(\beta)-\G_i^{(1)}))=O_p(1)$. By Lemma \ref{lemmaA1},
\begin{equation}
\label{bound-U5'}
\sup_{\|x\|=1}\sup_{\|y\|=1}\sup_{\beta \in B_n^*(r)}U_5'(\beta,\x,\y)\leq \pi_n \gamma_n^{(0)} O_p(1).
\end{equation}
Using \eqref{bound-U5}, \eqref{bound-U5'} and \eqref{bound-U}, we obtain by the hypotheses of the lemma that
$$\sup_{\|x\|=1}\sup_{\|y\|=1}\sup_{\beta \in B_n^*(r)}|U_5(\beta,\x,\y)|\leq  \{\big(\pi_n \gamma_n^{(0)}o_p(n)\big)^{1/2}+\big(\pi_n \gamma_n^{(0)}n E(\tau_n)\big)^{1/2}\} O_p(1)=o_p(1).$$

We now treat $\cE_n^{(2)}(\beta)$. Using \eqref{diagonal}, we see that
\begin{align*}
\x^T \H_n^{-1/2}\cE_n^{(2)}(\beta)\H_n^{-1/2}\y&=\sum_{i=1}^{n}\x^T \H_n^{-1/2} \X_i^T \F_i(\beta) \Ri^{-1} {\rm diag}(\X_i \H_n^{-1/2} \y) \G_i^{(2)}(\beta)\e_i^*\\
&=U_2(\x,\y)+U_4(\beta,\x,\y)+U_6(\beta,\x,\y),
\end{align*}
where
\begin{align*}
U_2(\x,\y)&=\sum_{i=1}^{n}\x^T \H_n^{-1/2} \X_i^T \F_i  \Ri^{-1} {\rm diag}(\X_i \H_n^{-1/2} \y) \G_i^{(2)}\e_i^*\\
U_4(\beta,\x,\y)&=\sum_{i=1}^{n}\x^T \H_n^{-1/2} \X_i^T (\F_i(\beta)-\F_i) \Ri^{-1} {\rm diag}(\X_i \H_n^{-1/2} \y) \G_i^{(2)}(\beta)\e_i^*\\
U_6(\beta,\x,\y)&= \sum_{i=1}^{n}\x^T \H_n^{-1/2} \X_i^T \F_i \Ri^{-1} {\rm diag}(\X_i \H_n^{-1/2} \y) (\G_i^{(2)}(\beta)-\G_i^{(2)})\e_i^*.
\end{align*}
By the Cauchy-Schwarz inequality \eqref{CS},
$|U_2(\x,\y)| \leq S_3(\beta_0,\x,\y)^{1/2}U^{1/2}$,
where $S_3(\beta,\x,\y)$ is given by \eqref{def-S3} and $U$ is the same as above. Using \eqref{S3-beta-xy} and \eqref{bound-U}, it follows that
$$\sup_{\|x\|=1}\sup_{\|y\|=1}|U_2(\x,\y)|=o_p(1).$$
Similarly, $|U_4(\beta,\x,\y)| \leq U_4'(\beta_0,\x,\y)^{1/2}U^{1/2}$,
where
\begin{align*}
U_4'(\beta,\x,\y)&=\sum_{i=1}^{n}\x^T \H_n^{-1/2}\X_i^T (\F_i(\beta)-\F_i) \Ri^{-1} {\rm diag}(\X_i \H_n^{-1/2}\y) \G_i^{(2)}(\beta)(\As_i)^{1/2}\Ri\\
& \quad \quad (\As_i)^{1/2}\G_i^{(2)}(\beta){\rm diag}(\X_i \H_n^{-1/2}\y)\Ri^{-1}(\F_i(\beta)-\F_i)\X_i \H_n^{-1/2} \x
\end{align*}
Using inequalities \eqref{ineq} and \eqref{diag-Xi}, it follows that
$$U_4'(\beta,\x,\y)\leq  \pi_n \gamma_n^{(0)}\max_{i\leq n}\lM^2(\G_i^{(2)}(\beta) (\As_i)^{1/2})T_1(\beta,\x,\x),$$
where $T_1(\beta,\x,\y)$ is given by \eqref{def-T1}. The matrix $\G_i^{(2)}(\beta) (\As_i)^{1/2}$ has $j$-th element given by:
$$\sqrt{\sigma_{i,jj}^*(\beta_0)}\,g_{ij}^{(2)}(\beta)=\frac{\sigma_{i,jj}^*(\beta_0)}{\sigma_{i,jj}^*(\beta)}
\cdot\left(\sqrt{\sigma_{i,jj}^*(\beta)}\,
g_{ij}^{(2)}(\beta)\right).$$
By Lemmas \ref{lemmaC} and \ref{lemmaD} (Appendix A),
\begin{equation}
\label{G2-As}
\sup_{\beta \in B_n^*(r)}\max_{i\leq n}\lM^2(\G_i^{(2)}(\beta) (\As_i)^{1/2})=O_p(1).
\end{equation}
 Using \eqref{bound-T1} and the fact that $\gamma_n^*=o_p(1)$, it follows that
$$\sup_{\beta \in B_n^*(r)} \sup_{\|\x\|=1} \sup_{\|\y\|=1}U_4'(\beta,\x,\y) \leq \pi_n^2 \gamma_n^{(0)} \gamma_n^* O_p(1)=\pi_n^2 \gamma_n^{(0)} o_p(1).$$
Using \eqref{bound-U} and the hypotheses of the lemma it follows that
$$\sup_{\beta \in B_n^*(r)} \sup_{\|\x\|=1} \sup_{\|\y\|=1}|U_4(\beta,\x,\y)|=o_p(1).$$

It remains to treat $U_6(\beta,\x,\y)$. By Cauchy-Schwarz inequality \eqref{CS},
$$|U_6(\beta,\x,\y)|\leq U_6'(\beta,\x,\y)^{1/2}U^{1/2},$$
where $U$ is the same as above and
\begin{align*}
U_6'(\beta,\x,\y)&=\sum_{i=1}^{n}\x^T \H_n^{-1/2} \X_i^T \F_i \Ri^{-1} {\rm diag}(\X_i \H_n^{-1/2} \y) (\G_i^{(2)}(\beta)-\G_i^{(2)})(\As_i)^{1/2}\Ri\\
& \quad \quad \quad (\As_i)^{1/2}(\G_i^{(2)}(\beta)-\G_i^{(2)}){\rm diag}(\X_i \H_n^{-1/2} \y)\Ri^{-1}\F_i \X_i \H_n^{-1/2}\x.
\end{align*}
Using inequalities \eqref{ineq} and \eqref{diag-Xi}, it follows that $U_6'(\beta,\x,\y)$ is less than or equal to
$$\pi_n \gamma_n^{(0)}\max_{i\leq n}\lM^2(\G_i^{(2)}(\A_i^*)^{1/2})\max_{i\leq n}\lM^2(\G_i^{(2)}(\beta)(\G_i^{(2)})^{-1}-\I)\cdot \x^T \H_n^{-1/2}\Hs_n \H_n^{-1/2}\x.$$
By \eqref{G2-As}, $\max_{i\leq n}\lM^2(\G_i^{(2)}(\A_i^*)^{1/2})=O_p(1)$. By Lemma \ref{lemmaA1}, $\sup_{x}\x^T \H_n^{-1/2}\Hs_n \H_n^{-1/2}\x=O_p(1)$. In can be shown that $\max_{i\leq n}\max_{i\leq n}\lM^2(\G_i^{(2)}(\beta)(\G_i^{(2)})^{-1}-\I)=O_p(1)$. Arguing as above, we infer that
$$\sup_{\beta \in B_n^*(r)} \sup_{\|\x\|=1} \sup_{\|\y\|=1}|U_6(\beta,\x,\y)|=o_p(1).$$
$\Box$

\section{Real-life example}
\label{section-ex}

In this section, we discuss an application of our method to a subset of the real-life dataset taken from \cite{sommer-katz-tarwotjo84}. This subset consists of $n=250$ preschool age rural Indonesian children which were examined every 3 months for 18 months for the presence of a respiratory disease. So each child was observed on $m=6$ occasions. The response $Y$ is a binary variable which takes value $1$ if the respiratory disease is present and value $0$ if the disease is absent. We consider a marginal logistic regression model (see Example \ref{example3}) with intercept parameter $\beta_0$ and 3 covariates: $X^{(1)}$ is a binary variable with values $0$ and $1$ giving the gender, $X^{(2)}$ is another binary variable with values 0 and 1 giving the vitamin A deficiency, and $X^{(3)}$ is the child's age in years at the beginning of the study, with possible values $1,2,\ldots,7$. In this case, $\mu(x)=e^x/(1+e^x)$. The model is
$${\rm logit}(Y_{ij})=\beta_0+\beta_1 X_{ij}^{(1)}+\beta_2 X_{ij}^{(2)}+\beta_3 X_{ij}^{(3)}, \quad i=1,\ldots,n, \quad j=1,\ldots,m.$$
Here $\beta=(\beta_0,\beta_1,\beta_2,\beta_3)$, so $p=4$.
We let $\X_{ij}=(X_{ij}^{(0)},X_{ij}^{(1)},X_{ij}^{(2)},X_{ij}^{(3)})$ where $X_{ij}^{(0)}=1$.

Since this data does not contain missing values,
we generated missingness indicator variables $I_{ij}$ using a Bernoulli distribution with probability $0.95$ of success. This gives $3.33\%$ missing responses.


We fit a logistic regression model with parameter $\gamma$ to the complete data set consisting of $(\I_i,\X_i)$ for $i=1,\ldots,250$, and we solved equation \eqref{eq50}.
The root of this equation is $\widehat{\gamma}=(3.514,0.025,0.391,-0.076)$. The estimates $\widehat{\pi}_{ij}$ for the missingness probabilities $\pi_{ij}$ are calcualted using the formula $\widehat{\pi}_{ij}=\pi_{ij}(\widehat{\gamma})$, where $\pi_{ij}(\gamma)$ is given by \eqref{def-pi-ij}. We compute the inverse probability weighted responses $$Y_{ij}^*=\frac{Y_{ij}I_{ij}}{\widehat{\pi_{ij}}} \quad i=1,\ldots,n, \ j=1,\ldots,m$$
and we solve the working independence GEE with weighted responses, which in this case is a system of 4 equations:
$$\sum_{i=1}^{n}\sum_{j=1}^{m}\X_{ij}^{(l)}\left(Y_{ij}^* -\frac{\exp(\X_{ij}^T \beta)}{1+\exp(\X_{ij}^T \beta)} \right)=0, \quad l=0,1,2,3$$
yielding the root $\beta^{\rm indep}=(-0.444,-0.552,0.258,-0.066)$.


We now apply our method to the dataset described above.
We start by computing the standardized values
$\widehat{Y}_{ij}=\widetilde{Y}_{ij}(\beta^{\rm indep})$, where
$$\widetilde{Y}_{ij}(\beta)=\frac{Y_{ij}^*-\mu_{ij}(\beta)}{\sqrt{\sigma_{i,jj}^*(\beta)} },$$
and $\sigma_{i,jj}^*(\beta)$ was calculated using \eqref{sigma_star_ijj} with $\pi_{ij}$ replaced by $\widehat{\pi}_{ij}$.

Recall that the conditional correlation matrix $\Rs_i$ of $\Y_i$ given $\X_i$ has elements:
$$r_{i,jk}^*=\frac{E[(Y_{ij}^*-\mu_{ij}(\beta))(Y_{ik}^*-\mu_{ik}(\beta))|\X_i]}{
\sqrt{\sigma_{i,jj}^*(\beta)} \cdot \sqrt{\sigma_{i,kk}^*(\beta)}}=E[\widetilde{Y}_{ij}(\beta) \widetilde{Y}_{ik}(\beta)|\X_i].$$
To estimate the matrix $\Rs_i$, we use the same matrix $\R_{i}(\alpha)=\widehat{\R}=(\widehat{r}_{jk})_{j,k=1,\ldots,m}$ for all $i$, with $\widehat{r}_{jj}=1$ for all $j=1,\ldots,m$, and for $j\not=k$, $\widehat{r}_{jk}$ are as in Examples 2 and 3 of \cite{liang-zeger86}:

Case 1: (1-dependent) $\widehat{r}_{jk}=0$ if $|j-k|\geq 2$ and $\widehat{r}_{j,j+1}=\widehat{r}_{j+1,j}=\widehat{\alpha}_j$
where
$$\widehat{\alpha}_{j}=\frac{1}{n-p} \sum_{i=1}^{n} \widehat{Y}_{ij} \widehat{Y}_{i,j+1}, \quad j=1,\ldots,m-1 $$
This produces the values $\widehat{\alpha}_1= 0.534, \widehat{\alpha}_2=0.559,\widehat{\alpha}_3=0.562,\widehat{\alpha}_4=0.443,\widehat{\alpha}_5=0.521$.

Case 2: (exhangeable) $\widehat{r}_{jk}=\widehat{\alpha}$ for all $j,k=1, \ldots,m$ with $j\not=k$, where
$$\widehat{\alpha}=\frac{1}{N-p}\sum_{i=1}^{n}\sum_{k=2}^{m}\sum_{j=1}^{k-1} \widehat{Y}_{ij} \widehat{Y}_{ik} \quad \mbox{with} \quad N=n \frac{m(m-1)}{2}$$
This produces the value $\widehat{\alpha}=0.492$.

Using these two cases, we solve equation \eqref{eq7}, taking into account that in the calculation of $\sigma_{i,jj}^{*}(\beta)$, $\pi_{ij}$ is replaced by $\widehat{\pi}_{ij}$. This consists of a system of 4 equations:
$$\sum_{i=1}^{n}\sum_{j=1}^{m}X_{ij}^{(l)} \frac{\mu'(\X_{ij}^T \beta)}{\sqrt{\sigma_{i,jj}^*(\beta)}}\,w_{jk} \, \frac{Y_{ik}^*-\mu(\X_{ik}^T \beta)}{\sqrt{\sigma_{i,kk}^*(\beta)}}=0, \quad l=0,1,2,3,$$
where ${\bf W}=\widehat{\R}^{-1}=(w_{jk})_{j,k=1,\ldots,m}$. We obtained the following estimates:
\begin{align*}
\widehat{\beta}^{(1)}&=(-0.448,-0.487,0.243,-0.068) \quad \mbox{for the $1$-dependent case}\\
\widehat{\beta}^{(2)}&=(-0.447,-0.550,0.256,-0.065) \quad \mbox{for the exchangeable case}.
\end{align*}

To evaluate the precision of these estimates, we compute the standard error of these estimates and the $p$-value of the two-sided test for $\beta=0$, using the asymptotic normality of $\widehat{\beta}$ given by Corollary \ref{asy-normality}:
\begin{equation}
\label{asy-norm-B}
{\bf B}^{-1/2}(\widehat{\beta}-\beta) \approx N_p({\bf 0},{\bf I}),
\end{equation}
where ${\bf B}=\H_n^{-1}\M_n \H_n^{-1}$. We estimate the matrix $\B$ by $\widehat{\bf B}=\widehat{\H}_n^{-1}\widehat{\M}_n \widehat{\H}_n^{-1}$, with matrices $\widehat{\M}_n$ and $\widehat{\H}_n$ computed as in Remark \ref{practice-remark}.
From \eqref{asy-norm-B}, we deduce that $\widehat{\beta}-\beta$ has approximately a $p$-variate normal distribution with mean vector ${\bf 0}$ and covariance matrix $\widehat{{\bf B}}$. Hence, for $l=0,1,2,3$, $\widehat{\beta}^{(l)}-\beta_l \approx N(0,b_l)$, where $b_{l}$ is the $l$-th element on the diagonal of $\widehat{{\bf B}}$.
It follows that the standard error (s.e.) of $\widehat{\beta}^{(l)}$ is $s\{\widehat{\beta}^{(l)}\}=\sqrt{b_l}$ and the $p$-value of the test of $H_0:\beta_l=0$ versus $H_1:\beta_l\not=0$ is $2P(Z>|\widehat{\beta}^{(l)}/\sqrt{b_l} |)$. In Table 1, we report the estimates, their standard errors and
$p$-values for the two examples of correlation matrices considered above ($1$-dependent and exchangeable).
\begin{center}
\begin{table}
\begin{tabular}{|c|c|c|c|||c|c|c|c|c} \hline
{\em $1$-dependent}                   & estimate & s.e. & $p$-value & {\em exchangeable} & estimate & s.e. & $p$-value  \\ \hline
intercept  & -0.448 & 0.268 & 0.095 & intercept & -0.447 & 0.272 & 0.101 \\
gender  & -0.487 & 0.219 & 0.026 & gender & -0.550 & 0.220 & 0.012  \\
vitamin A  & 0.243 & 0.224 & 0.279 & vitamin A & 0.256 &  0.226 & 0.256 \\
age    & -0.068  & 0.056 & 0.220 & age  & -0.065 & 0.056 & 0.241 \\
\hline
\end{tabular}
\caption{Parameter estimates, standard errors (s.e.) and $p$-values for Case 1 ($1$-dependent) and Case 2 (exchangeable) correlation matrices}
\end{table}
\end{center}

We conclude that at a 5\% significance level, we reject the hypothesis $\beta_1=0$, but we do not have enough evidence to reject the hypothesis $\beta_2=0$ or the hypothesis $\beta_3=0$. This means that the gender seems to have a significant effect on the presence of respiratory disease, but vitamin A deficiency and age do not influence the presence of this disease.

\appendix

\section{Some auxiliary results}

In this appendix, we gather some auxiliary results which were used in the proofs of Lemmas \ref{lemmaA1}, \ref{lemmaA2} and \ref{lemmaA3}.

\begin{lemma}
\label{lemmaB}
Suppose Assumption (AH) holds. If $\gamma_n^* \stackrel{P}{\to}0$, then
\begin{equation}
\label{mu-prime}
\sup_{\beta \in B_n^*(r)}\max_{i\leq n}\max_{j \leq m}\left|\frac{\mu'(\X_{ij}^T \beta)}{\mu'(\X_{ij}^T \beta_0)}-1 \right|=O_{p}((\gamma_n^*)^{1/2})
\end{equation}
\begin{equation}
\label{mu-square}
\sup_{\beta \in B_n^*(r)}\max_{i\leq n}\max_{j \leq m}\left|\frac{\mu^2(\X_{ij}^T \beta)}{\mu^2(\X_{ij}^T \beta_0)}-1 \right|=O_{p}((\gamma_n^*)^{1/2}).
\end{equation}
\end{lemma}

\noindent {\bf Proof:} This follows by Taylor's formula,
using the fact that $k_{n}^{(1)}=O_p(1)$, respectively $k_{n}^{(0)}=O_p(1)$. See also Lemma B.1 of \cite{XY}. $\Box$

\vspace{3mm}

\begin{lemma}
\label{lemmaC}
Suppose Assumption (AH) holds. If $\gamma_n^* \stackrel{P}{\to}0$ then
$$\sup_{\beta \in B_n^*(r)} \max_{i \leq n} \max_{j \leq m} \left|\frac{\sigma_{i,jj}^*(\beta)}{\sigma_{i,jj}^*(\beta_0)}-1\right|=O_{p}((\gamma_n^{*})^{1/2}).$$
\end{lemma}

\noindent {\bf Proof:} Note that $\frac{\sigma_{i,jj}^*(\beta)}{\sigma_{i,jj}^*(\beta_0)}-1$ is equal to
$$\frac{\mu'(\X_{ij}^T \beta)-\mu'(\X_{ij}^T \beta_0)+\left(\frac{1}{\pi_{ij}}-1 \right)\big(\mu'(\X_{ij}^T \beta)-\mu'(\X_{ij}^T \beta_0)+\mu^2(\X_{ij}^T\beta)-\mu^2(\X_{ij}^T\beta_0)\big)}{\mu'(\X_{ij}^T \beta_0)+\left(\frac{1}{\pi_{ij}}-1 \right)\big(\mu'(\X_{ij}^T \beta_0)+\mu^2(\X_{ij}^T\beta_0)\big)}.$$
Since $\mu'$ is non-negative and $\pi_{ij} \leq 1$,
$$\left|\frac{\sigma_{i,jj}^*(\beta)}{\sigma_{i,jj}^*(\beta_0)}-1\right| \leq
2\left|\frac{\mu'(\X_{ij}^T \beta)-\mu'(\X_{ij}^T \beta_0)}{\mu'(\X_{ij}^T \beta_0)}\right|+\left|\frac{\mu^2(\X_{ij}^T \beta)-\mu^2(\X_{ij}^T\beta_0)}{\mu^2(\X_{ij}^T\beta_0)}\right|.$$
The conclusion follows by Lemma \ref{lemmaB}. $\Box$

\begin{lemma}
\label{lemmaD}
Suppose Assumptions (AH) and (M) hold. Then
$$\sup_{\beta \in B_n^*(r)} \max_{i \leq n} \max_{j \leq m} \left(\sqrt{\sigma_{i,jj}^*(\beta)}\,|g_{ij}^{(2)}(\beta)|\right)=O_p(1).$$
\end{lemma}

\noindent {\bf Proof:} Recalling definition \eqref{def-gij-2} of $g_{ij}^{(2)}(\beta)$, we see that:
\begin{equation}
\label{sigma-g2}
\sqrt{\sigma_{i,jj}^*(\beta)}\,g_{ij}^{(2)}(\beta)=
-\frac{2\left(\frac{1}{\pi_{ij}}-1\right)\mu(\X_{ij}^T \beta)\mu'(\X_{ij}^T \beta)+\frac{1}{\pi_{ij}} \mu''(\X_{ij}^T\beta)}{2\sigma_{i,jj}^*(\beta)},
\end{equation}
and hence
$$\sqrt{\sigma_{i,jj}^*(\beta)}\,|g_{ij}^{(2)}(\beta)|\leq \left( \frac{1}{\pi_{ij}}-1 \right)\frac{|\mu(\X_{ij}^T \beta)\mu'(\X_{ij}^T \beta)|}{\sigma_{i,jj}^*(\beta)}+\frac{1}{\pi_{ij}}\cdot \frac{|\mu''(\X_{ij}^T\beta)|}{2\sigma_{i,jj}^*(\beta)}.$$
In the first term on the right-hand side of this inequality, we use the fact that
$\sigma_{i,jj}^*(\beta) \geq
\big(\frac{1}{\pi_{ij}}-1\big)(\mu'(\X_{ij}^T \beta)+\mu^2(\X_{ij}^T \beta)) \geq \big(\frac{1}{\pi_{ij}}-1\big)\mu^2(\X_{ij}^T \beta)$, whereas for the second term we use the fact that $\sigma_{i,jj}^*(\beta)\geq \mu'(\X_{ij}^T \beta)$. We obtain that
\begin{align*}
\sqrt{\sigma_{i,jj}^*(\beta)}\,|g_{ij}^{(2)}(\beta)| & \leq \left|\frac{\mu'(\X_{ij}^T\beta)}{\mu(\X_{ij}^T \beta)}\right|+\frac{1}{2\pi_{ij}}\left|\frac{\mu''(\X_{ij}^T\beta)}{\mu'(\X_{ij}^T \beta)} \right|\leq k_n^{(0)}+\frac{1}{2}\rho_n k_n^{(1)}.
\end{align*}
The conclusion follows by Assumptions (AH) and (M). $\Box$.

\begin{lemma}
\label{lemmaF}
Suppose Assumptions (AH) and (M) hold.
If $\gamma_n^{*} \stackrel{P}{\to}0$ then
$$\sup_{\beta \in B_n^*(r)} \max_{i \leq n} \max_{j \leq m} \left|\frac{f_{ij}(\beta)}{f_{ij}(\beta_0)}-1\right|=O_{p}((\gamma_n^{*})^{1/2}).$$
\end{lemma}

\noindent {\bf Proof:}  Recall that $\frac{\partial}{\partial \beta^T}f_{ij}(\beta)=g_{ij}^{(1)}(\beta)\X_{ij}^T$. By Taylor's formula, for any $\beta \in B_n^{*}(r)$, there exists $\beta_{ij} \in B_n^*(r)$ such that
$f_{ij}(\beta)-f_{ij}(\beta_0)=g_{ij}^{(1)}(\beta_{ij})\X_{ij}^T(\beta-\beta_0)$. Since
$$\|\X_{ij}^T(\beta-\beta_0)\| \leq \|\X_{ij}^T \H_n^{-1/2}\| \cdot \|\H_n^{1/2}(\beta-\beta_0)\|\leq (\gamma_n^{(0)})^{1/2}(\tau_n)^{1/2}r=(\gamma_n^*)^{1/2} r,$$
it follows that
\begin{equation}
\label{fij-minus-1}
\left|\frac{f_{ij}(\beta)}{f_{ij}(\beta_0)}-1 \right|\leq  \left| \frac{g_{ij}^{(1)}(\beta_{ij})}{f_{ij}(\beta_{ij})}\right| \cdot \left| \frac{f_{ij}(\beta_{ij})}{f_{ij}(\beta_0)}\right| (\gamma_n^*)^{1/2} r.
\end{equation}

\noindent By definition \eqref{def-fij} of $f_{ij}(\beta)$, we have
$\frac{f_{ij}(\beta)}{f_{ij}(\beta_0)}=
\frac{\mu'(\X_{ij}^T \beta)}{\mu'(\X_{ij}^T \beta_0)}\cdot \frac{\sqrt{\sigma_{i,jj}^*(\beta_0)}}{\sqrt{\sigma_{i,jj}^*(\beta)}}$.
By Lemmas \ref{lemmaB} and \ref{lemmaC},
\begin{equation}
\label{fij-Op}
\sup_{\beta \in B_n^*(r)}\max_{i\leq n}\max_{j\leq m} \left| \frac{f_{ij}(\beta)}{f_{ij}(\beta_0)}\right|=O_{p}(1).
\end{equation}

\noindent A direct calculation based on definition \eqref{def-gij-1} of $g_{ij}^{(1)}(\beta)$ and relation \eqref{sigma-g2} shows that
$$\frac{g_{ij}^{(1)}(\beta)}{f_{ij}(\beta)}= \frac{\mu''(\X_{ij}^T \beta)}{\mu'(\X_{ij}^T \beta)}+\sqrt{\sigma_{i,jj}^*(\beta)}\,g_{ij}^{(2)}(\beta).$$
By Assumption (AH) and Lemma \ref{lemmaD}, it follows that
\begin{equation}
\label{gij1-fij}
\sup_{\beta \in B_n^*(r)}\max_{i\leq n}\max_{j \leq m}\left|\frac{g_{ij}^{(1)}(\beta)}{f_{ij}(\beta)} \right|=O_{p}(1).
\end{equation}
The conclusion follows from relations \eqref{fij-minus-1}, \eqref{fij-Op} and \eqref{gij1-fij}. $\Box$


\end{document}